\documentclass[11pt]{article}
\usepackage{amssymb,amsmath,url,fullpage,rotating} 
\usepackage{authblk}

\def\eqdef{:=}
\def\Ber{\mathrm{Ber}}

\def\CS{\mathrm{CS}}
\def\KL{\mathrm{KL}}
\def\JS{\mathrm{JS}}
\def\dmu{\mathrm{d}\mu}
\def\bbR{\mathbb{R}}

\def\calX{\mathcal{X}}
\def\calF{\mathcal{F}}
\def\dmu{\mathrm{d}\mu}
\def\Bhat{\mathrm{Bhat}}

\def\tr{\mathrm{tr}}
\def\calE{\mathcal{E}}
\def\calN{\mathcal{N}}
\def\alphasum{\alpha_{\sum}}

\newtheorem{lemma}{Lemma}
\newtheorem{theorem}{Theorem}
\newtheorem{example}{Example}

\title{Cumulant-free closed-form formulas for some common (dis)similarities between densities of an exponential family}

\date{}

\author[$\star$]{Frank Nielsen}
\affil[$\star$]{Sony Computer Science Laboratories, Inc.}
\affil[$\star$]{Tokyo, Japan}
\affil[$\star$]{{\small\tt Frank.Nielsen@acm.org}}
\author[$\dagger$]{Richard Nock}
\affil[$\dagger$]{CSIRO Data61 \& Australian National University}
\affil[$\dagger$]{Sydney, Australia}
\affil[$\dagger$]{{\small\tt Richard.Nock@data61.csiro.au}}

\begin{document}
\maketitle
\begin{abstract}
It is well-known that the Bhattacharyya, Hellinger, Kullback-Leibler, $\alpha$-divergences, and Jeffreys' divergences between densities belonging to a same exponential family have generic closed-form formulas relying on the strictly convex and real-analytic cumulant function characterizing the exponential family. 
In this work, we report (dis)similarity formulas which bypass the explicit use of the cumulant function and highlight the role of 
quasi-arithmetic means and their multivariate mean operator extensions. 
In practice, these cumulant-free formulas are handy when implementing these (dis)similarities using legacy Application Programming Interfaces (APIs) since our method requires only to partially factorize the densities canonically of the considered exponential family.
\end{abstract}


\section{Introduction}

Let $(\calX,\calF,\mu)$ be a measure space~\cite{PM-1995} with sample space $\calX$, $\sigma$-algebra of events $\calF$, and positive measure $\mu$ (i.e., Lebesgue or counting measures).
The Kullback-Leibler divergence~\cite{KL-1951}  (KLD), Jeffreys' divergence~\cite{Jeffreys-1946} (JD), Bhattacharyya coefficient (BC), Bhattacharyya distance~\cite{Bhattacharyya-1946,Kailath-1967} (BD) and Hellinger distance~\cite{Hellinger-1909} (HD) between 
two probability measures $P$ and $Q$ dominated by $\mu$ with respective Radon-Nikodym densities  $p=\frac{\mathrm{d}P}{\dmu}$ and
 $q=\frac{\mathrm{d}Q}{\dmu}$ 
 are statistical (dis)similarities defined respectively by:
\begin{align}
D_\KL[p:q] &:= \int p(x)\log\frac{p(x)}{q(x)}\dmu(x),&\mbox{(Kullback-Leibler divergence)}\\
D_J[p,q] &:= D_\KL[p:q]+D_\KL[q:p],&\mbox{(Jeffreys' divergence)}\\
\rho[p,q] &:= \int_\calX \sqrt{p(x)\ q(x)} \dmu(x),&\mbox{(Bhattacharyya coefficient)}\\
D_B[p,q] &:= -\log(\rho[p,q]),&\mbox{(Bhattacharyya distance)}\\
D_H[p,q] &:= \sqrt{1-\rho[p,q]}. &\mbox{(Hellinger distance)}
\end{align}

KLD is an oriented distance (i.e., $D_\KL[p:q] \not=D_\KL[q:p]$). 
JD is a common symmetrization of the KLD which is not a metric distance because JD fails the triangle inequality.
BC is a similarity which ensures that $\rho[p,q]\in(0,1]$.
BD is a symmetric non-metric distance, and HD is a metric distance.
(KLD and JD are homogeneous divergences.)

All these divergences require to calculate definite integrals which can be calculated in theory using 
Risch pseudo-algorithm~\cite{risch1970solution} which depends on an oracle.
In practice, computer algebra systems like {\tt Maple}\textregistered{} or {\tt Maxima} implement different subsets of 
 the theoretical Risch pseudo-algorithm and thus face different limitations.

When the densities $p$ and $q$ belong to a same parametric family $\calE=\{p_\theta\ :\ \theta\in\Theta\}$ of densities, 
i.e., $p=p_{\theta_1}$ and $q=p_{\theta_2}$, these (dis)similarities yield {\em equivalent} parameter (dis)similarities.
For example, we get the parameter divergence $P(\theta_1:\theta_2)=D_\KL^\calE(\theta_1:\theta_2):=D_\KL(p_{\theta_1}:p_{\theta_2})$.
We use notationally the brackets to indicate that the (dis)similarities parameters are densities, 
and the parenthesis to indicate parameter (dis)similarities.

In particular, when $\calE$ is a natural exponential family~\cite{EF-2009,Barndorff-2014} (NEF) 
\begin{equation}
\calE= \left\{p_\theta(x)=\exp\left(\theta^\top t(x)-F(\theta)+k(x)\right) \ : \ \theta\in\Theta\right\},  \label{eq:factorization}
\end{equation}
with $t(x)$ denoting the sufficient statistics, $k(x)$ an auxiliary measure carrier term, and 
\begin{equation}\label{eq:cumulantF}
F(\theta) := \log\left(\int_\calX \exp(\theta^\top t(x)) \dmu(x)\right),
\end{equation}
the cumulant function\footnote{More precisely, $\kappa_\theta(u)=F(\theta+u)-F(\theta)$ is the cumulant generating function of the sufficient statistic $t(x)$ from which all moments can be recovered (the cumulant generating function being the logarithm of the moment generating function).} also called  log-normalizer, log-partition function, free energy, or log-Laplace transform.
Parameter $\theta$ is called the natural parameter.
The cumulant function is a strictly smooth convex function (real analytic function~\cite{Barndorff-2014}) defined on the open convex natural parameter space $\Theta$.
Let $D$ denote the dimension of the parameter space $\Theta$ (i.e., the order of the exponential family) and $d$ the dimension of the sample space $\calX$.
We further assume full {\em regular} exponential families~\cite{Barndorff-2014} so that $\Theta$ is a 
non-empty open convex domainand $t(x)$ is a minimal sufficient statistic~\cite{nosteepsinglytruncnormal-1994}.

Many common families of distributions $\{p_\lambda(x)\ \lambda\in\Lambda\}$ are exponential families in disguise after reparameterization: 
$p_\lambda(x)=p_{\theta(\lambda)}(x)$. 
Those families are called exponential families (i.e., EFs and not natural EFs to emphasize that $\theta(u)\not=u$), and their densities are canonically factorized as follows:
\begin{equation}\label{eq:canef}
p(x;\lambda)=\exp\left(\theta(\lambda)^\top t(x)-F(\theta(\lambda))+k(x)\right).
\end{equation}
We  call parameter $\lambda$ the  source parameter (or the ordinary parameter, with $\lambda\in\Lambda$, the source/ordinary parameter space) and parameter $\theta(\lambda)\in\Theta$ is the corresponding natural parameter.
Notice that the canonical parameterization of Eq.~\ref{eq:canef} is {\em not} unique: 
For example, adding a constant term $c\in\bbR$  to $F(\theta)$ can be compensated by subtracting this constant to $k(x)$, or multiplying the sufficient statistic  $t(x)$ by a symmetric invertible matrix $A$ can be compensated by multiplying $\theta(\lambda)$ by the inverse of $A$ so that $(A\theta(\lambda)))^\top (A^{-1}t(x))=\theta(\lambda))^\top AA^{-1}t(x) =  \theta(\lambda)^\top t(x)$.

Another useful parameterization of exponential families is the {\em moment parameter}~\cite{EF-2009,Barndorff-2014}: 
$\eta=E_{p_\theta}[t(x)]=\nabla F(\theta)$. The moment parameter space shall be denoted by $H$:
\begin{equation}
H=\left\{ E_{p_\theta}[t(x)]\ :\ \theta\in\Theta\right\}.
\end{equation}
Let $C=\overline{\mathrm{CH}}\{ t(x)\ :\ x\in\calX\}$ be the closed convex hull of the range of the sufficient statistic.
When $H=\mathrm{int}(C)$, the family is said {\em steep}~\cite{nosteepsinglytruncnormal-1994}.
In the remainder, we consider full regular and steep exponential families.

To give one example, consider the family of univariate normal densities:
\begin{equation}
\mathcal{N}:=\left\{ p(x;\lambda)=\frac{1}{\sigma \sqrt{2 \pi}} \exp\left(-\frac{1}{2}\left(\frac{x-\mu}{\sigma}\right)^{2}\right),\ \lambda=(\mu,\sigma^2)\in\bbR\times\bbR_{++}\right\}.
\end{equation}
Family $\calN$ is interpreted as an exponential family of order $D=2$ with univariate parametric densities ($d=1$), indexed by the source \ parameter $\lambda=(\mu,\sigma^2)\in\Lambda$ with $\Lambda=\bbR\times\bbR_{++}$.
The corresponding natural parameter is $\theta(\lambda)=\left(\frac{\mu}{\sigma^2},-\frac{1}{2\sigma^2}\right)$,
 and the moment parameter is $\eta(\lambda)=(E_{p_\lambda}[x],E_{p_\lambda}[x^2])=(\mu,\mu^2+\sigma^2)$ since $t(x)=(x,x^2)$.
The cumulant function for the normal family is $F(\theta)=-\frac{\theta_{1}^{2}}{4 \theta_{2}}+\frac{1}{2} \log \left(-\frac{\pi}{\theta_{2}}\right)$.

When densities $p=p_{\theta_1}$ and $q=p_{\theta_2}$ both belong to the same exponential family, 
we have the following well-known closed-form expressions~\cite{BR-2011,EF-2009} for  the  (dis)similarities introduced formerly:
\begin{align}
D_\KL[p_{\theta_1}:p_{\theta_2}] &= {B_F}^*(\theta_1:\theta_2), &\mbox{Kullback-Leibler divergence}\\
D_J[p_{\theta_1}:p_{\theta_2}] &= (\theta_2-\theta_1)^\top (\eta_2-\eta_1),&\mbox{Jeffreys' divergence}\\
\rho[p_{\theta_1}:p_{\theta_2}] &= \exp(-J_F(\theta_1:\theta_2)),&\mbox{Bhattacharyya coefficient}\\
D_B[p_{\theta_1}:p_{\theta_2}] &= J_F(\theta_1:\theta_2)&\mbox{Bhattacharyya distance}\\
D_H[p_{\theta_1}:p_{\theta_2}] &= \sqrt{1-\exp(-J_F(\theta_1:\theta_2))},&\mbox{Hellinger distance}
\end{align}
where $D^*$ indicates the reverse divergence $D^*(\theta_1:\theta_2):=D(\theta_2:\theta_1)$ (parameter swapping\footnote{Here, the star `*' is not to be confused with the Legendre-Fenchel transform used in information geometry~\cite{IG-2016}. 
For a Bregman divergence $B_F(\theta_1:\theta_2)$, we have the reverse Bregman divergence
 ${B_F}^*(\theta_1:\theta_2)=B_{F}(\theta_2:\theta_1)=B_{F^*}(\nabla F(\theta_1):\nabla F(\theta_2))$, where $F^*$ denotes convex conjugate of $F$ obtained by the Legendre-Fenchel transform.}), 
and $B_F$ and $J_F$ are the Bregman divergence~\cite{Bregman-1967} and the Jensen divergence~\cite{BR-2011} induced by the functional generator $F$, respectively:
\begin{eqnarray}
B_F(\theta_1:\theta_2) &:=& F(\theta_1)-F(\theta_2)-(\theta_1-\theta_2)^\top \nabla F(\theta_2)\\
J_F(\theta_1:\theta_2) &:=& \frac{F(\theta_1)+F(\theta_2)}{2} - F\left(\frac{\theta_1+\theta_2}{2}\right).
\end{eqnarray}

In information geometry, the Bregman divergences are the canonical divergences of dually flat spaces~\cite{nielsen2018elementary}.

More generally, the Bhattacharrya distance/coefficient can be skewed with a parameter $\alpha\in(0,1)$ to yield
the $\alpha$-skewed Bhattacharyya distance/coefficient~\cite{BR-2011}:
\begin{eqnarray}
D_{B,\alpha}[p:q] &:=& -\log(\rho_\alpha[p:q]),\\
\rho_\alpha[p:q] &:=& \int_\calX p(x)^{\alpha}  q(x)^{1-\alpha}\dmu(x),\\ \label{eq:BCalpha}
&=&-\log \int_\calX q(x) \left(\frac{p(x)}{q(x)}\right)^{\alpha}\dmu(x).\label{eq:bhatalphafdiv}
\end{eqnarray}
The ordinary Bhattacharyya distance and coefficient are recovered for $\alpha=\frac{1}{2}$.
In statistics, the maximum skewed Bhattacharrya distance is called Chernoff information~\cite{Chernoff-1952,Chernoff-2013}   used in Bayesian hypothesis testing:
\begin{equation}
D_C[p:q]:= \max_{\alpha\in(0,1) } D_{B,\alpha}[p:q].
\end{equation}

Notice that the Bhattacharyya skewed $\alpha$-coefficient of Eq.~\ref{eq:BCalpha} also appears in the definition of the $\alpha$-divergences~\cite{IG-2016} (with $\alpha\in\bbR$):
\begin{equation}
D_\alpha[p:q] := \left\{
\begin{array}{ll}
\frac{1}{\alpha(1-\alpha)}  \left(1-\rho_\alpha[p:q]\right), & \alpha\in\bbR\backslash\{0,1\} \label{eq:alphanormalized}\\
D_1[p:q]=D_\KL[p:q], & \alpha=1\\
D_0[p:q]=D_\KL[q:p],& \alpha=0
\end{array}
\right.
\end{equation}
The  $\alpha$-divergences belong to the class of $f$-divergences~\cite{Csiszar-1967} which are the invariant divergences\footnote{A statistical invariant divergence $D$ is such that $D[p_{\theta}:p_{\theta+\mathrm{d}\theta}]=\sum_{i,j} g_{ij}(\theta)\mathrm{d}\theta_i\mathrm{d}\theta_j$, where $I(\theta)=[ g_{ij}(\theta)]$ is the Fisher information matrix~\cite{IG-2016}.} in information geometry~\cite{IG-2016}.

When densities $p=p_{\theta_1}$ and $q=p_{\theta_2}$ both belong to the same exponential family $\calE$, we get 
the following closed-form formula~\cite{BR-2011}:
\begin{eqnarray}
D_{B,\alpha}[p_{\theta_1}:p_{\theta_2}] &=& J_{F,\alpha}(\theta_1:\theta_2), 
\end{eqnarray}
where $J_{F,\alpha}$ denotes the  $\alpha$-skewed Jensen divergence:
\begin{equation}
J_{F,\alpha}(\theta_1:\theta_2) := \alpha F(\theta_1)+(1-\alpha) F(\theta_2)- F( \alpha\theta_1+(1-\alpha)\theta_2).
\end{equation}
All these closed-form formula can be obtained from the calculation of the following generic integral~\cite{SM-2011}:
\begin{equation} \label{eq:integralI}
I_{\alpha,\beta}[p:q] := \int p(x)^{\alpha}q(x)^\beta\dmu(x),
\end{equation}
when  $p=p_{\theta_1}$ and $q=p_{\theta_2}$.
Indeed, provided that $\alpha\theta_1+\beta\theta_2\in\Theta$, we have
\begin{equation} 
I_{\alpha,\beta}[p_{\theta_1}:p_{\theta_2}] =
\exp\left(- (\alpha F(\theta_1)+\beta F(\theta_2)-F(\alpha\theta_1+\beta\theta_2)) \right)E_{p_{\alpha\theta_1+\beta\theta_2}}\left[\exp((\alpha+\beta-1)k(x))\right].\label{eq:solintegralI}
\end{equation}
The calculation of $I_{\alpha,\beta}$ in Eq.~\ref{eq:solintegralI} is easily achieved by bypassing the computation of the antiderivative of the integrand in Eq.~\ref{eq:integralI} (using the fact that $\int p_\theta(x)\dmu(x)=1$ for any $\theta\in\Theta$), see~\cite{SM-2011}.

In particular, we get the following special cases:
\begin{itemize}
	\item When $\alpha+\beta=1$, $I_{\alpha,1-\alpha}[p_{\theta_1}:p_{\theta_2}]=\exp(-J_{F,\alpha}(\theta_1:\theta_2))$ since $\alpha\theta_1+(1-\alpha)\theta_2\in\Theta$ (domain $\Theta$ is convex).
	
	\item When $k(x)=0$, and $\alpha+\beta>0$, $I_{\alpha,\beta}[p_{\theta_1}:p_{\theta_2}]=\exp\left(- (\alpha F(\theta_1)+\beta F(\theta_2)-F(\alpha\theta_1+\beta\theta_2)) \right)$. 
	This is always the case when $\Theta$ is a convex cone (e.g., Gaussian or Wishart family), see~\cite{EFMixture-2012}.
	
	\item When $\alpha+\beta=1$ with arbitrary $\alpha$, 
	\begin{equation}
	I_{\alpha,1-\alpha}[p_{\theta_1}:p_{\theta_2}]=\exp\left(- (\alpha F(\theta_1)+(1-\alpha)F(\theta_2)-F(\alpha\theta_1+(1-\alpha)\theta_2)) \right).
	\end{equation}
	This setting is useful for getting truncated series of $f$-divergences when the exponential family has an affine space $\Theta$, see~\cite{FDiv-2013,Fdivpowerchi-2019}.
	
\end{itemize}

When $\alpha\rightarrow 1$ or $\alpha\rightarrow 0$, we get the following limits of the $\alpha$-skewed Bhattacharrya distances:
\begin{eqnarray}
\lim_{\alpha\rightarrow 1} \frac{1}{\alpha(1-\alpha)} D_{B,\alpha}[p:q] &=& D_\KL[p:q]\\
\lim_{\alpha\rightarrow 0} \frac{1}{\alpha(1-\alpha)} D_{B,\alpha}[p:q] &=& D_\KL[q:p].
\end{eqnarray}

It follows that when the densities $p=p_{\theta_1}$ and $q=p_{\theta_2}$ both belong to the same exponential family, we obtain 
\begin{eqnarray}
\lim_{\alpha\rightarrow 0} \frac{1}{\alpha(1-\alpha)} J_{F,\alpha}(\theta_1:\theta_2) &=& B_F(\theta_2:\theta_1),\\
\lim_{\alpha\rightarrow 1} \frac{1}{\alpha(1-\alpha)} J_{F,\alpha}(\theta_1:\theta_2) &=& B_F(\theta_1:\theta_2).
\end{eqnarray}

In practice, we would like to get closed-form formula for the (dis)similarities when the densities belong to the same exponential families using the source reparameterization $\lambda\in\Lambda$:
\begin{eqnarray}
D_\KL[p_{\lambda_1}:p_{\lambda_2}] &=& {B_F}(\theta(\lambda_2):\theta(\lambda_1)),\\
D_J[p_{\lambda_1}:p_{\lambda_2}] &=& (\theta(\lambda_1):\theta(\lambda_2))^\top (\eta(\lambda_2)-\eta(\lambda_1)),\\
\rho[p_{\lambda_1}:p_{\lambda_2}] &=& \exp(-J_F(\theta(\lambda_1):\theta(\lambda_2))),\\
D_B[p_{\lambda_1}:p_{\lambda_2}] &=& J_F(\theta(\lambda_1):\theta(\lambda_2))\\
D_H[p_{\lambda_1}:p_{\lambda_2}] &=& \sqrt{1-\rho(\theta(\lambda_1):\theta(\lambda_2))},
\end{eqnarray}
where $\theta(\cdot)$ and $\eta(\cdot)$ are the $D$-variate functions for converting the source parameter $\lambda$ to the natural parameter $\theta$ and the moment parameter $\eta$, respectively.
The Chernoff information between two densities $p_{\lambda_1}$ and $p_{\lambda_2}$ of the same exponential family amounts to a Jensen-Chernoff divergence~\cite{Chernoff-2013}:
\begin{eqnarray}
\mathrm{JC}({\lambda_1}:{\lambda_2}) := D_C[p_{\lambda_1}:p_{\lambda_2}] &=& \max_{\alpha\in(0,1) } D_{B,\alpha}[p_{\lambda_1}:p_{\lambda_2}],\\
&=&  \max_{\alpha\in(0,1) }  J_{F,\alpha}(\theta(\lambda_1):\theta(\lambda_2)),
\end{eqnarray}
that is the maximal value of a skew Jensen divergence for $\alpha\in(0,1)$.

Thus to have closed-form formula, we need to explicit both the $\theta(\cdot)$ conversion function and the cumulant function $F$.
This can be prone to human calculus mistakes (e.g., report manually these formula without calculus errors for the multivariate Gaussian family). Furthermore, the cumulant function $F$ may not be available in closed-form~\cite{PEF-2016}.

In this work, we show how to easily bypass the explicit use of the cumulant function $F$.
Our method is based on a simple trick, and makes the programming of these (dis)similarities easy using off-the-shelf functions of application programming interfaces (APIs) (e.g., the density function, the entropy function or the moment function of a distribution family).

This paper is organized as follows:
Section~\ref{sec:BC} explains the method for the Bhattacharyya coefficient and its related dissimilarities.
Section~\ref{sec:KL} further carry on the principle of bypassing the explicit use of the cumulant function and its gradient for the calculation of the Kullback-Leibler divergence and its related  Jeffreys' divergence.
Section~\ref{sec:concl} summarizes the results.
Throughout the paper, we present several examples to illustrate the methods.
Appendix~\ref{sec:CAS} displays some code written using the computer algebra system (CAS) {\sc Maxima} to recover some formula for some exponential families. Appendix~\ref{sec:furtherexamples} provides further examples.

\section{Cumulant-free formula for the Bhattacharyya coefficient and  distances derived thereof}\label{sec:BC}

\subsection{A method based on a simple trick}

The densities of an exponential family  have all the same support~\cite{Barndorff-2014} $\calX$. 
 Consider  {\em any} point $\omega\in\calX$ in the support. 
Then observe that we can write the cumulant of a natural exponential family as:
\begin{equation}
F(\theta)=-\log p_{\theta}(\omega)+t(\omega)^\top \theta +  k(\omega).
\end{equation}

Since the generator $F$ of a Jensen or Bregman divergence is defined modulo an affine function $a^\top \theta+b$ (i.e., $J_F=J_{G}$ and $B_F=B_G$ for $G(\theta)=F(\theta)+a^\top \theta+b$ for $a\in\bbR^D$ and $b\in\bbR$), 
we consider the following {\em equivalent} generator (term $+t(\omega)^\top \theta +  k(\omega)$ is affine) expressed using the density parameterized by $\lambda\in\Lambda$:
\begin{equation}
F_\lambda(\lambda) = F(\theta(\lambda)) = -\log(p_\lambda(\omega)).
\end{equation}

Then the Bhattacharyya coefficient is expressed by this {\em cumulant-free expression} using the source parameterization $\lambda$:

\begin{eqnarray}\label{eq:bhatcoeff}
\forall \omega\in\calX,\quad \rho[p_{\lambda_1},p_{\lambda_2}] =  
   \frac{\sqrt{p_{\lambda_1}(\omega)p_{\lambda_2}(\omega)}}{p_{\bar\lambda}(\omega)} = 
	\sqrt{\frac{p_{\lambda_1}(\omega)}{p_{\bar\lambda}(\omega)}} \sqrt{\frac{p_{\lambda_2}(\omega)}{p_{\bar\lambda}(\omega)}} 
,
\end{eqnarray}
with
\begin{equation}
\bar\lambda=\lambda\left(\frac{\theta(\lambda_1)+\theta(\lambda_2)}{2}\right).
\end{equation}

Similarly, the Bhattacharyya  distance is written as
\begin{eqnarray}\label{eq:bhat}
\forall \omega\in\calX,\quad D_B[p_{\lambda_1},p_{\lambda_2}] =  
 \log\left( \frac{p\left(\omega;\bar\lambda\right)}{\sqrt{p(\omega,\lambda_1)p(\omega,\lambda_2)}}\right).
\end{eqnarray}

Let $l(x;\lambda):=\log p(x;\lambda)$ be the log-density.
Then we have
\begin{eqnarray}\label{eq:logbhat}
\forall \omega\in\calX,\quad D_B[p_{\lambda_1},p_{\lambda_2}] =  
 l\left(\omega;\bar\lambda\right) - \frac{l(\omega,\lambda_1)+l(\omega,\lambda_2)}{2}.
\end{eqnarray}
This is a Jensen divergence for the strictly convex function $-l(x;\theta)$ (wrt. $\theta$) since $-l(x;\theta)\equiv F(\theta)$ (modulo  an affine term).

Thus we do not need to {\em explicitly} use the expression of the cumulant function $F$ in Eq.~\ref{eq:bhatcoeff} and Eq.~\ref{eq:bhat} but  we need the following {\em parameter $\lambda\Leftrightarrow\theta$ conversion functions}: 
\begin{enumerate}
\item $\theta(\lambda)$ the ordinary-to-natural parameter conversion function, and 
\item its {\em reciprocal function} $\lambda(\theta)$ (i.e., $\lambda(\cdot)=\theta^{-1}(\cdot)$), the natural-to-source parameter conversion function
\end{enumerate}
so that we can calculate the ordinary $\lambda$-parameter $\bar\lambda$ corresponding to the natural mid-parameter $\frac{\theta(\lambda_1)+\theta(\lambda_2)}{2}$:
\begin{equation}
\bar\lambda=\lambda\left(\frac{\theta(\lambda_1)+\theta(\lambda_2)}{2}\right).
\end{equation}

Notice that in general, a linear interpolation in the natural parameter $\theta$ corresponds to a {\em non-linear interpolation} in the source  parameterization $\lambda$ when $\theta(u)\not=u$.

Since $\lambda(\cdot)=\theta^{-1}(\cdot)$, we can interpret the non-linear interpolation $\bar\lambda$ as a generalization of {\em quasi-arithmetic mean}~\cite{BR-2011} (QAM):
\begin{equation}
\bar\lambda = \theta^{-1}\left(\frac{\theta(\lambda_1)+\theta(\lambda_2)}{2}\right)=: M_{\theta}(\lambda_1,\lambda_2),
\end{equation}
where 
\begin{equation}\label{eq:qam}
M_f(a,b):=f^{-1}\left(\frac{1}{2}f(a)+\frac{1}{2}f(b)\right),
\end{equation}
is the quasi-arithmetic mean induced by a strictly monotonic and smooth function $f$.\footnote{Notice that 
$M_f(a,b)=M_g(a,b)$ when $g(u)=c_1 f(u)+ c_2$ with $c_1\not=0\in\bbR$ and $c_2\in\bbR$.
Thus we can assume that $f$ is a strictly increasing and smooth function. Function $f$ is said to be in standard form if $f(1)=0$, $f'(u)>0$ and $f'(1)=1$).}

We can extend the quasi-arithmetic mean to a weighted quasi-arithmetic mean as follows:
\begin{equation}\label{eq:qawm}
M_{f,\alpha}(a,b):=f^{-1}\left(\alpha f(a)+(1-\alpha)f(b)\right),\quad \alpha\in [0,1].
\end{equation}
Weighted quasi-arithmetic means are strictly monotone~\cite{WOpDiss-2014} when the $\mathrm{Range}(f)\subset\bbR$. 

Let us remark that extensions of weighted quasi-arithmetic means have been studied recently in information geometry to describe geodesic paths~\cite{PathQAM-2011,PathQAM-2015} $\gamma_{ab}=\{ M_{f,\alpha}(a,b)\ :\ \alpha\in[0,1]\}$.

In 1D, a bijective function on an interval (e.g., a parameter conversion function in our setting)  is a strictly monotonic function, and thus $f$ defines a quasi-arithmetic mean.

To define quasi-arithmetic mean with multivariate generators using Eq.~\ref{eq:qam}, we need to properly define $f^{-1}$.
When the function $f$ is separable, i.e., $f(x)=\sum_{i=1}^D f_i(x_i)$ with the $f_i$'s strictly monotone and smooth functions, we can  straightforwardly define the multivariate mean 
as $M_f(x,x'):=(M_{f_1}(x_1,x_1'),\ldots, M_{f_D}(x_1,x_1'))$.

In general, the inverse function theorem in multivariable calculus~\cite{clarke1976inverse} states 
that there exists {\em locally} an inverse function $f^{-1}$ for a continuously differentiable function $f$ at point $x$ if
the determinant of the Jacobian $J_f(x)$ of $f$ at $x$ is non-zero. 
Moreover, $f^{-1}$ is continuously differentiable and we have $J_{f^{-1}}(y)=[J_f(x)]^{-1}$ (matrix inverse) at $y=f(x)$.

In some cases, we get the existence of a global inverse function.
For example, when $f=\nabla H$ is the gradient of Legendre-type strictly convex and smooth function $H$, the reciprocal function $f^{-1}$
is well-defined $f^{-1}=\nabla H^*$ and global, where $H^*$ denotes the convex conjugate of $H$ (Legendre-type).
In that case, $f=\nabla H$ is a strictly monotone operator\footnote{A mapping $M:\bbR^d\rightarrow\bbR^d$ is a strictly monotone operator iff $(p-q)^\top (M(p)-M(q))> 0$ for all $p,q\in\bbR^d$ with $p\not=q$. 
Monotone operators are a generalization of the concept of univariate monotonic functions.} since it is the gradient of a strictly convex function.
We also refer to~\cite{OperatorMeans-2011} for some work on operator means, and to~\cite{QAM-Covar-2009} for
multivariate quasi-arithmetic means of covariance matrices.

To summarize, we can compute the Bhattacharyya coefficient (and Bhattacharyya/Hellinger distances) using the parametric density function 
$p_\lambda(x)$ and a quasi-arithmetic mean $\bar\lambda=M_{\theta}(\lambda_1,\lambda_2)$ on the source parameters $\lambda_1$ and 
$\lambda_2$ as:

\begin{eqnarray}\label{eq:bhatcoeffcf}
\forall\omega\in\calX,\quad
\rho[p_{\lambda_1},p_{\lambda_2}] =  
   \frac{\sqrt{p_{\lambda_1}(\omega)p_{\lambda_2}(\omega)}}{p(\omega;M_\theta(\lambda_1,\lambda_2))}
\end{eqnarray}
 using the notation $p(x;\theta):=p_\theta(x)$, and
\begin{equation}\label{eq:bhatcf}
\forall\omega\in\calX,\quad
D_B[p_{\lambda_1},p_{\lambda_2}] =  
 \log\left( \frac{p\left(\omega;M_\theta(\lambda_1,\lambda_2)\right)}{\sqrt{p(\omega,\lambda_1)p(\omega,\lambda_2)}}\right).
\end{equation}

Similarly, we get the following cumulant-free expression for the Hellinger distance:

\begin{equation}
\forall\omega\in\calX,\quad
D_H[p_{\lambda_1},p_{\lambda_2}]
= \sqrt{1-
\frac{\sqrt{p_{\lambda_1}(\omega)p_{\lambda_2}(\omega)}}{p_(\omega;M_\theta(\lambda_1,\lambda_2))}
}.
\end{equation}

The Hellinger distance proves useful when using some generic algorithms which require to handle metric distances.
For example, the $2$-approximation factor of Gonzalez~\cite{Gonzalez-1985} for $k$-center metric clustering.

These cumulant-free formula are all the more convenient as in legacy software API, we usually have  access to the density function of the probability family.
Thus if a parametric family of an API is an exponential family $\calE$, 
we just need to implement the corresponding quasi-arithmetic mean $M_\theta^\calE$.

More generally, the $\alpha$-skewed Bhattacharyya distance~\cite{BR-2011} for $\alpha\in(0,1)$ is expressed using
the following cumulant-free expression:
\begin{equation}\label{eq:bhatalpha}
\forall\omega\in\calX,\quad
D_{B,\alpha}[p_{\lambda_1}:p_{\lambda_2}] =
\log\left( \frac{p(\omega,\lambda_1)^{\alpha} p(\omega,\lambda_2)^{1-\alpha}}{p\left(\omega;M_{\theta,\alpha}(\lambda_1,\lambda_2)\right)} \right)
\end{equation}
with 
\begin{equation}
M_{\theta,\alpha}(\lambda_1,\lambda_2) = \lambda(\alpha\theta(\lambda_1)+(1-\alpha)\theta(\lambda_2)) =: \bar\lambda_\alpha.
\end{equation}

Notice that the geometric $\alpha$-barycenter of two densities $p_{\theta_1}$ and $p_{\theta_2}$ of an exponential family $\calE$ is a scale density of $\calE$:
$$
p(x;\alpha\theta_1+(1-\alpha)\theta_2)= \left( \frac{p(\omega;\alpha\theta_1+(1-\alpha)\theta_2)}{p^\alpha(\omega;\theta_1)p^{1-\alpha}(\omega;\theta_2)} \right)  p^\alpha(x;\theta_1)p^{1-\alpha}(x;\theta_2),\quad \forall\omega\in\calX.
$$

Let $\tilde{p}_\theta(x)=\exp(t(x)^\top\theta+k(x))$ denote the unnormalized density of an exponential family (so that $p_\theta(x)=\tilde{p}_\theta(x)\exp(-F(\theta))$).
We have the following invariant:
\begin{equation}
 \forall\omega\in\calX,\quad
\frac{\tilde{p}(\omega;{M_\theta,\alpha}(\lambda_1,\lambda_2))}{\tilde{p}(\omega,\lambda_1)^{\alpha} \tilde{p}(\omega,\lambda_2)^{1-\alpha}} = 1.
\end{equation}

It follows that we have:
\begin{eqnarray}
\rho_\alpha[p_{\theta_1}:p_{\theta_2}] &=& \frac{{p}(\omega;\alpha\theta_1+(1-\alpha)\theta_2)}{p(\omega,\lambda_1)^{\alpha} p(\omega,\lambda_2)^{1-\alpha}},\\    
 &=& \frac{\exp(F(\alpha\theta_1+(1-\alpha)\theta_2))}{\exp(\alpha F(\theta_1)) \exp((1-\alpha) F(\theta_2))},\\
&=& \exp(-J_{F,\alpha}(\theta_1:\theta_2)).
\end{eqnarray}

The Cauchy-Schwarz divergence~\cite{jenssen2006cauchy,nielsen2017holder} is defined by
\begin{equation}
D_\CS(p,q) \eqdef -\log\left( \frac{\int_{\calX} p(x) q(x) \dmu(x)}{\sqrt{\left(\int_{\calX} p(x)^{2} \dmu(x)\right)\left(\int_{\calX} q(x)^{2}\dmu(x)\right)}} \right) \geq 0.
\end{equation}
Thus the Cauchy-Schwarz divergence is a {\em projective divergence}:   $D_\CS(p,q) =D_\CS(\lambda p,\lambda q)$ for any $\lambda>0$. 
It can be shown that the Cauchy-Schwarz divergence between two densities $p_{\theta_1}$ and $p_{\theta_2}$ of an exponential family with a natural parameter space a cone~\cite{OnicescuCSDiv-2020} (e.g., Gaussian) is:
\begin{equation}
D_\CS(p_{\theta_1},p_{\theta_2})=J_F(2\theta_1:2\theta_2)+\log\left( \frac{\sqrt{E_{p_{2\theta_1}}\left[\exp(k(x))\right]E_{p_{2\theta_2}}\left[\exp(k(x))\right]}}{E_{p_{\theta_1+\theta_2}}\left[\exp(k(x))\right]}\right).
\end{equation}
See~\cite{OnicescuCSDiv-2020} for the formula extended to H\"older divergences which generalize the Cauchy-Schwarz divergence.

\subsection{Some illustrating examples}
 
Let us start with an example of a continuous exponential family which relies on the arithmetic mean $A(a,b)=\frac{a+b}{2}$:

\begin{example}\label{ex:expofamily}
Consider the family of exponential distributions with rate parameter $\lambda>0$.
The densities of this continuous EF writes as $p_\lambda(x)=\lambda\exp(-\lambda x)$ with  support $\calX=[0,\infty)$.
From the partial canonical factorization  of densities following Eq.~\ref{eq:factorization}, we get that
 $\theta(u)=u$ and $\theta^{-1}(u)=u$ so that $M_\theta(\lambda_1,\lambda_2)=A(\lambda_1,\lambda_2)=\frac{\lambda_1+\lambda_2}{2}$ is the arithmetic mean.
Choose $\omega=0$ so that $p_\lambda(\omega)=\lambda$.
It follows that
\begin{eqnarray*}
\rho[p_{\lambda_1},p_{\lambda_2}] &=& \frac{\sqrt{p(\omega,\lambda_1)p(\omega,\lambda_2)}}{p\left(\omega;\bar\lambda\right)},\\
&=& \frac{\sqrt{\lambda_1\lambda_2}}{\frac{\lambda_1+\lambda_2}{2}},\\
&=& \frac{2\sqrt{\lambda_1\lambda_2}}{\lambda_1+\lambda_2}.
\end{eqnarray*}

Let $G(a,b)=\sqrt{ab}$ denote the geometric mean for $a,b>0$.
Notice that since $G(\lambda_1,\lambda_2)\leq A(\lambda_1,\lambda_2)$, we have $\rho[p_{\sigma_1},p_{\sigma_2}] \in (0,1]$.
The Bhattacharyya distance between two exponential densities of rate $\lambda_1$ and $\lambda_2$ is
$$
D_\Bhat[p_{\sigma_1},p_{\sigma_2}] = -\log \rho[p_{\lambda_1},p_{\lambda_2}]= \log\frac{\lambda_1+\lambda_2}{2} -\frac{1}{2}\log\lambda_1\lambda_2\geq 0.
$$
Since the logarithm function is monotonous, we have $\log  A(\lambda_1,\lambda_2)\geq \log G(\lambda_1,\lambda_2)$ and therefore
 we check that $D_\Bhat[p_{\sigma_1},p_{\sigma_2}]\geq 0$.
\end{example}

Next, we consider a discrete exponential family which exhibits the geometric mean:

\begin{example}\label{ex:Poissonfamily}
The Poisson family of probability mass functions (PMFs) $p_\lambda(x)=\frac{\lambda^x\exp(-\lambda)}{x!}$ where $\lambda>0$ denotes the intensity parameter and $x\in\calX=\{0,1,\ldots, \}$ is a discrete exponential family 
with $t(x)=x$ ($\omega=0$, $t(\omega)=0$ and $p_\lambda(\omega)=\exp(-\lambda)$), $\theta(u)=\log u$ 
and $\lambda(u)=\theta^{-1}(u)=\exp(u)$. 
Thus the quasi-arithmetic mean associated with the Poisson family
 is $M_\theta(\lambda_1,\lambda_2)=G(\lambda_1,\lambda_2)=\sqrt{\lambda_1\lambda_2}$ the geometric mean.
It follows that the Bhattacharrya coefficient is
\begin{eqnarray*}
\rho[p_{\lambda_1},p_{\lambda_2}] &=& \frac{p\left(\omega;\bar\lambda\right)}{\sqrt{p(\omega,\lambda_1)p(\omega,\lambda_2)}},\\
&=& \frac{\exp(-\sqrt{\lambda_1\lambda_2})}{\sqrt{\exp(-\lambda_1) \exp(-\lambda_2)}} = \exp\left(\frac{\lambda_1+\lambda_2}{2}-\sqrt{\lambda_1\lambda_2}\right)
\end{eqnarray*}

Hence, we recover the Bhattacharrya distance between two Poisson pmfs:
$$
D_B[p_{\lambda_1},p_{\lambda_2}] = \frac{\lambda_1+\lambda_2}{2}-\sqrt{\lambda_1\lambda_2}\geq 0.
$$

The negative binomial distribution with known number of failures yields also the same natural parameter and geometric mean (but the density $p(\omega,\lambda)$ is different).
\end{example}

To illustrate the use of the power means $P_r(a,b)=(a^r+b^r)^{\frac{1}{r}}$ of order $r\in\bbR$ (also called H\"older means) for $r\not=0$ (with $\lim_{r\rightarrow 0} P_r(a,b)=G(a,b)$),  let us consider the family of Weibull distributions.

\begin{example}\label{ex:Weibullfamily}
The Weibull distributions with a prescribed shape parameter $k \in \bbR_{++}$ (e.g., exponential distributions when $k=1$, Rayleigh distributions when $k=2$) form an exponential family.
The density of a Weibull distribution with scale $\lambda$ and fixed shape parameter $k$ is
$$
p_\lambda(x)= \frac{k}{\lambda}\left(\frac{x}{\lambda}\right)^{k-1} e^{-(x / \lambda)^{k}},\quad x\in\calX=\bbR_{++}
$$

The ordinary$\leftrightarrow$natural parameter conversion functions are
$\theta(u)=\frac{1}{u^k}$ and $\theta^{-1}(u)=\frac{1}{u^{\frac{1}{k}}}$.
It follows that 
$\bar\lambda=M_\theta(\lambda_1,\lambda_2)=P_{-k}(\lambda_1,\lambda_2)= \left(\frac{1}{2}\lambda_1^{-k}+\frac{1}{2}\lambda_2^{-k}\right)^{-\frac{1}{k}}$ 
is the power means of order $-k$.
We choose $\omega=1$ and get $p_\lambda(\omega)=\frac{k}{\lambda^k} e^{-\frac{1}{\lambda^{k}}}$.
It follows the closed-form formula for the Bhattacharrya coefficient for integer $k\in\{2,\ldots , \}$:
\begin{eqnarray}
\rho[p_{\lambda_1},p_{\lambda_2}] &=& {\sqrt{p(\omega,\lambda_1)p(\omega,\lambda_2)}}\frac{p\left(\omega;\bar\lambda\right)},\\
&=& \frac{2\sqrt{\lambda_1^k\lambda_2^k}}{\lambda_1^k+\lambda_2^k}.
\end{eqnarray}
For $k=2$, the Weibull family yields the Rayleigh family with $\lambda=\sqrt{2}\sigma$.
The Bhattacharyya coefficient is $\rho[p_{\lambda_1},p_{\lambda_2}]=\frac{G(\lambda_1,\lambda_2)}{Q(\lambda_1,\lambda_2)}$
 where $G$ and $Q$ denotes the geometric mean and the quadratic mean, respectively (with $Q\geq G$).
\end{example}

Table~\ref{tab:case1d} summarizes the quasi-arithmetic means  associated with common univariate exponential families.
Notice that a same quasi-arithmetic mean can be associated to many exponential families:
For example, the Gaussian family with fixed variance or the exponential distribution family have both $\theta(u)=u$ yielding the arithmetic mean.

\begin{table}
\centering
$$
\begin{array}{llll}
\mbox{Exponential family $\calE$} & \theta(u) & \theta^{-1}(u) & \bar\lambda=M_{\theta}(\lambda_1,\lambda_2)\\ \hline
\mbox{Exponential} &  u & u & A(\lambda_1,\lambda_2)=\frac{\lambda_1+\lambda_2}{2} \\
\mbox{Poisson} &  \log u & \exp(u) & G(\lambda_1,\lambda_2)=\sqrt{\lambda_1\lambda_2} \\
\mbox{Laplace (fixed $\mu$)} & \frac{1}{u} & \frac{1}{u} &  H(\lambda_1,\lambda_2)=\frac{2\lambda_1\lambda_2}{\lambda_1+\lambda_2}\\
\mbox{Weibull (fixed shape $k>0$)} & \frac{1}{u^k} & \frac{1}{u^{\frac{1}{k}}} & 
P_{\frac{1}{k}}(\lambda_1,\lambda_2)=(\frac{1}{2}\lambda_1^{\frac{1}{k}}+ \frac{1}{2}\lambda_1^{\frac{1}{k}})^{-\frac{1}{k}} \\
\mbox{Bernoulli} &  \log\frac{u}{1-u} & \frac{\exp(u)}{1+\exp(u)} & \Ber(\lambda_1,\lambda_2)=\frac{a_{12}}{1+a_{12}}\\
 & \mbox{(logit function)} & \mbox{(logistic function)}  & \mbox{with\ } a_{12}=\sqrt{\frac{\lambda_1\lambda_2}{(1-\lambda_1)(1-\lambda_2)}}\\
\hline
\end{array}
$$
\caption{Quasi-arithmetic means $M_\theta(\lambda_1,\lambda_2)=\theta^{-1}\left(\frac{\theta(\lambda_1)+\theta(\lambda_2)}{2}\right)$ associated with some common discrete and continuous exponential families of order $D=1$.\label{tab:case1d}}
\end{table}

Let us now consider multi-order exponential families.
We start with the bi-order exponential family  of Gamma distributions.

\begin{example}\label{ex:Gamma}
The density of a Gamma distribution is
$$
p(x;\alpha,\beta)=\frac{\beta^{\alpha}}{\Gamma(\alpha)} x^{\alpha-1} e^{-\beta x},\quad x\in\calX=(0,\infty),
$$
for a shape parameter $\alpha>0$ and rate parameter $\beta>0$ (i.e., a $2$-order exponential family with $\lambda=(\lambda_1,\lambda_2)=(\alpha,\beta)$).
The natural parameter is $\theta(\lambda)=(\lambda_1-1,-\lambda_2)$ and the inverse function
is $\lambda(\theta)=(\theta_1+1,-\theta_2)$.
It follows that the generalized quasi-arithmetic mean is the bivariate arithmetic mean: 
$M_\theta(\lambda_1,\lambda_2)=A(\lambda_1,\lambda_2)=(\frac{\alpha_1+\alpha_2}{2},\frac{\beta_1+\beta_2}{2})=(\bar\alpha,\bar\beta)$.
We choose $\omega=1$ so that $p(\omega;\alpha,\beta)=\frac{\beta^{\alpha}}{\Gamma(\alpha)}  e^{-\beta}$.

We get the Bhattacharrya coefficient:
\begin{eqnarray*}
\rho[p_{\alpha_1,\beta_1}:p_{\alpha_2,\beta_2}] &=& \frac{\sqrt{p(\omega,\lambda_1)p(\omega,\lambda_2)}}{p\left(\omega;\bar\lambda\right)},\\
&=& \sqrt{\frac{\beta_1^{\alpha_1}\beta_2^{\alpha_2}}{\Gamma{\alpha_1}\Gamma{\alpha_2}}} \bar\beta^{-\bar\alpha},
\end{eqnarray*}
and the Bhattacharrya distance:
$$
D_\alpha[p_{\alpha_1,\beta_1}:p_{\alpha_2,\beta_2}] = \bar\alpha\log\bar\beta-\frac{\alpha_1\log\beta_1+\alpha_2\log\beta_2}{2}
+\log \frac{\sqrt{\Gamma(\alpha_1)\Gamma(\alpha_2)}}{\Gamma(\bar\alpha)}.
$$
\end{example}

The Dirichlet family which exhibits a separable (quasi-)arithmetic mean:

\begin{example}\label{ex:DirichletBhat}
Consider the family of Dirichlet distributions with densities defined on the $(d-1)$-dimensional open standard simplex  support 
\begin{equation}
\calX=\Delta_d:=\left\{(x_1,\ldots, x_d)\ :\ x_{i} \in(0,1),\ \sum_{i=1}^d x_{i}=1\right\}.
\end{equation}
The family of Dirichlet distributions including the family of Beta distributions when $d=2$.
The density of a Dirichlet distribution is defined  by:
\begin{equation}
p_\alpha(x)= \frac{\Gamma\left(\sum_{i=1}^d \alpha_{i}\right)}{\prod_{i=1}^d 
\Gamma\left(\alpha_{i}\right)} \prod_{i=1}^d x_{i}^{\alpha_{i}-1}.
\end{equation}

The Dirichlet distributions  and are used in in Bayesian statistics as the conjugate priors  of the multinomial family.
The Dirichlet distributions form an exponential family with $d$-dimensional natural parameter $\theta=(\alpha_1-1,\ldots,\alpha_d-1)$ ($D=d$) and vector of
sufficient statistics $t(x)=(\log x_1,\ldots, \log x_d)$.
The induced quasi-arithmetic mean is a multivariate separable arithmetic means, i.e., the multivariate arithmetic mean 
$\bar\lambda=M_\theta(\alpha_1,\alpha_2)=A(\alpha_1,\alpha_2)=\frac{\alpha_1+\alpha_2}{2}$.

Let us choose $\omega=\left(\frac{1}{d},\ldots, \frac{1}{d}\right)$ so that 
$p(\omega;\alpha)= \frac{\Gamma\left(\sum_{i=1}^d \alpha_{i}\right)}{\prod_{i=1}^d 
\Gamma\left(\alpha_{i}\right)}  \frac{1}{d^{(\sum_{i=1}^d \alpha_{i})-d}}$.

We get the following Bhattacharrya coefficient between two Dirichlet densities $p_{\alpha_1}$ and $p_{\alpha_2}$:
\begin{equation}
\rho[p_{\alpha_1}:p_{\alpha_2}] = 
\frac{B\left(\frac{\alpha_1+\alpha_2}{2}\right)}{\sqrt{B(\alpha_1) B(\alpha_2)}},
\end{equation}
where 
\begin{equation}
B(\alpha):=\frac{\prod_{i=1}^d \Gamma\left(\alpha_{i}\right)}{\Gamma\left(\sum_{i=1}^d \alpha_{i}\right)}.
\end{equation}

It follows that the Bhattacharrya distance  between two Dirichlet densities $p_{\alpha_1}$ and $p_{\alpha_2}$ is 
\begin{equation}
D_B[p_{\alpha_1}:p_{\alpha_2}] = \log\left(
\frac{\sqrt{B(\alpha_1) B(\alpha_2)}}{B\left(\frac{\alpha_1+\alpha_2}{2}\right)}
\right).
\end{equation}
The $\alpha$-skewed Bhattacharrya coefficient for a scalar $\alpha\in(0,1)$ is:
\begin{equation}
\rho_\alpha[p_{\alpha_1}:p_{\alpha_2}] = 
\frac{B\left(\alpha\alpha_1+(1-\alpha)\alpha_2\right)}{B^\alpha(\alpha_1) B^{1-\alpha}(\alpha_2)},
\end{equation}
and the $\alpha$-skewed Bhattacharrya distance:
\begin{eqnarray}
D_{B,\alpha}[p_{\alpha_1}:p_{\alpha_2}] &=& \log\left(
\frac{B^\alpha(\alpha_1) B^{1-\alpha}(\alpha_2)}{B\left(\alpha\alpha_1+(1-\alpha)\alpha_2\right)}
\right),\\
&=& \alpha \log B(\alpha_1)+(1-\alpha)\log B(\alpha_2)-\log B\left(\alpha\alpha_1+(1-\alpha)\alpha_2\right),
\end{eqnarray}
with 
\begin{equation}
\log B(\alpha)=\sum_{i=1}^d \log\Gamma(\alpha_i)-\log\Gamma\left(\sum_{i=1}^d \alpha_i\right).
\end{equation}
(This is in accordance with Eq. 15--17 of~\cite{Dirichlet-2008}.)
\end{example}

Finally, we consider the case of the multivariate Gaussian family:

\begin{example}\label{ex:GaussianBhat}
Consider the $d$-dimensional Gaussian family~\cite{JS-2019} also called MultiVariate Normal (MVN) family.
The parameter $\lambda=(\lambda_v,\lambda_M)$ of a MVN consists of a vector part $\lambda_v=\mu$ and a $d\times d$ matrix part $\lambda_M=\Sigma$.
The Gaussian density is given by
$$
p_\lambda(x;\lambda) =  \frac{1}{(2\pi)^{\frac{d}{2}}\sqrt{|\lambda_M|}}  \exp\left(-\frac{1}{2} (x-\lambda_v)^\top \lambda_M^{-1} (x-\lambda_v)\right),
$$ 
where $|\cdot|$ denotes the matrix determinant.

Partially factorizing the density into the canonical form of exponential family, we find that
$\theta(\lambda)=\left(\lambda_M^{-1}\lambda_v,\frac{1}{2}\lambda_M^{-1}\right)$ and 
$\lambda(\theta)=(\frac{1}{2}\theta_M^{-1}\theta_v,\frac{1}{2}\theta_M^{-1})$.
It follows that the multivariate weighted mean is $M_\theta^\alpha(\lambda_1,\lambda_2)=\bar\lambda_\alpha=(\mu_\alpha,\Sigma_\alpha)$ with  
\begin{equation}
\Sigma_\alpha=  \left(\alpha\Sigma_1^{-1}+(1-\alpha) \Sigma_2^{-1}\right)^{-1},
\end{equation}
the matrix harmonic barycenter and 
\begin{equation}
\mu_\alpha= \Sigma_\alpha \left(\alpha\Sigma_1^{-1}\mu_1+(1-\alpha) \Sigma_2^{-1}\mu_2\right).
\end{equation}

We choose $\omega=0$ with
$p_\lambda(0;\lambda) =  \frac{1}{(2\pi)^{\frac{d}{2}}\sqrt{|\lambda_M|}}  \exp\left(-\frac{1}{2} \lambda_v^\top \lambda_M^{-1}\lambda_v\right)$.
Let $\Delta_\mu=\mu_2-\mu_1$. 
It follows the following closed-form formula for the Bhattacharyya coefficient between Gaussian densities:

$$
\rho[p_{\mu_1,\Sigma_1},p_{\mu_2,\Sigma_2}] = 
\left|\frac{\Sigma_1+\Sigma_2}{2}\right|^{-\frac{1}{2}} |\Sigma_1|^{\frac{1}{4}} |\Sigma_2|^{\frac{1}{4}} 
\exp\left(-\frac{1}{8} \Delta_\mu^\top  \left(\frac{\Sigma_1+\Sigma_2}{2}\right)^{-1} \Delta_\mu\right).  
$$

Thus the Bhattacharrya distance is
$$
D_{B,\alpha}[p_{\mu_1,\Sigma_1},p_{\mu_2,\Sigma_2}] = 
\frac{1}{2}\left(\alpha\mu_1^\top\Sigma_1^{-1}\mu_1+(1-\alpha)\mu_2^\top\Sigma_2^{-1}\mu_2-\mu_\alpha^\top\Sigma_\alpha^{-1}\mu_\alpha
+\log \frac{|\Sigma_1|^{\alpha}|\Sigma_2|^{1-\alpha}}{|\Sigma_\alpha|}\right),
$$
and the Hellinger distance:
$$
D_H[p_{\mu_1,\Sigma_1},p_{\mu_2,\Sigma_2}] = 
\sqrt{1-\frac{|\Sigma_{1}|^{\frac{1}{4}} 
|\Sigma_{2}|^{\frac{1}{4}}}{\left|\frac{\Sigma_{1}+\Sigma_{2}}{2}\right|^{\frac{1}{2}1}} 
\exp \left(-\frac{1}{8}\Delta_\mu\top \left|\frac{\Sigma_{1}+\Sigma_{2}}{2}\right|^{-1} \Delta_\mu\right) }.
$$

The Cauchy-Schwarz divergence between two Gaussians~\cite{OnicescuCSDiv-2020} is: 
\begin{eqnarray}
D_\CS(p_{\mu_1,\Sigma_1},p_{\mu_2,\Sigma_2})  &=& \frac{1}{2}\log\left(\frac{1}{2^d} \frac{\sqrt{|\Sigma_1|\ \Sigma_2||}}{|(\Sigma_1^{-1}+\Sigma_2^{-1})^{-1}|} \right)\nonumber\\
&& +\frac{1}{2}\mu_1^\top\Sigma_1^{-1}\mu_1+\frac{1}{2}\mu_2^\top\Sigma_2^{-1}\mu_2\nonumber\\
&& -\frac{1}{2}(\Sigma_1^{-1}\mu_1+\Sigma_2^{-1}\mu_2)^\top (\Sigma_1^{-1}+\Sigma_2^{-1})^{-1}  (\Sigma_1^{-1}\mu_1+\Sigma_2^{-1}\mu_2).
\end{eqnarray}

\end{example}

\section{Cumulant-free formula for the Kullback-Leibler divergence and related divergences}\label{sec:KL}

The Kullback-Leibler divergence (KLD) $D_\KL[p:q] :=\int p(x)\log\frac{p(x)}{q(x)}\dmu(x)$ between two densities $p$ and $q$ also called  the relative entropy~\cite{CoverThomasIT-2012} amounts to a reverse Bregman divergence, $D_\KL[p_{\theta_1}:p_{\theta_2}] = {B_F}^*(\theta_1:\theta_2)=B_F(\theta_2:\theta_1)$, when the densities belong to the same exponential family $\calE$, where the Bregman generator $F$ 
is the cumulant function of $\calE$.

We present below two techniques to calculate the KLD by avoiding to compute the integral: 
\begin{itemize}
\item The first technique, described in \S\ref{sec:KLlimit}, considers the KLD as a limit case of $\alpha$ skewed Bhattacharrya distance. 
\item The second technique relies on the availability of off-the-shelf formula for the entropy and moment of the sufficient statistic (\S\ref{sec:KLhmoment}), and is derived using the Legendre-Fenchel divergence.
\end{itemize}

\subsection{Kullback-Leibler divergence as the limit case of a skewed Bhattacharrya distance}\label{sec:KLlimit}

We can obtain closed-form formula for the Kullback-Leibler divergence by considering the limit case of $\alpha$ skewed Bhattacharrya distance:
\begin{eqnarray}
D_\KL[p_{\lambda_1}:p_{\lambda_2}] &=& \lim_{\alpha\rightarrow 0} \frac{1}{\alpha} B_\alpha[p_{\lambda_2}:p_{\lambda_1}],\\
&=& \lim_{\alpha\rightarrow 0} \frac{1}{\alpha} \log\rho_\alpha[p_{\lambda_2}:p_{\lambda_1}],\\
&=&  \lim_{\alpha\rightarrow 0} \frac{1}{\alpha} \log \frac{p\left(\omega;M_{\theta,\alpha}(\lambda_2,\lambda_1)\right)}{p(\omega,\lambda_2)^{\alpha} p(\omega,\lambda_1)^{1-\alpha}},\\
&=&  \log\left(\frac{p(\omega;\lambda_1)}{p(\omega;\lambda_2)}\right) 
+ \lim_{\alpha\rightarrow 0}   \frac{1}{\alpha} \log \left(\frac{p\left(\omega;M_{\theta,\alpha}(\lambda_2,\lambda_1)\right)}{p(\omega,\lambda_1)}\right).\label{eq:klefalpha}
\end{eqnarray}

When we deal with uni-order exponential families ($D=1$),
we can use a first-order Taylor expansion of the quasi-arithmetic means when $\alpha\simeq 0$ (see~\cite{alphadiv-2020}):
\begin{equation}
M_\theta(a,b) =_{\alpha\simeq 0} 
a +\alpha \frac{\theta(b)-\theta(a)}{\theta'(a)}
+ o\left(\alpha(\theta(b)-\theta(a))\right),
\end{equation}
where $\theta'(\cdot)$ denote the derivative of the ordinary-to-natural parameter conversion function.

It follows that we have:
\begin{equation}
D_\KL[p_{\lambda_1}:p_{\lambda_2}] = \log\left(\frac{p(\omega;\lambda_1)}{p(\omega;\lambda_2)}\right)
+ \lim_{\alpha\rightarrow 0} \frac{1}{\alpha} 
\log\left(\frac{
p\left(\omega;\lambda_1+\alpha \frac{\theta(\lambda_2)-\theta(\lambda_1)
}{\theta'(\lambda_1)}\right)}{p(\omega;\lambda_1)}\right), \forall \omega\in\calX \label{eq:klefalpha2}
\end{equation}

Notice that we need to calculate {\em case by case} the limit as it depends on the density expression $p(x;\lambda)$ of the exponential family. This limit can be computed symbolically using a computer algebra system (e.g., using {\sc Maxima}\footnote{\url{http://maxima.sourceforge.net/}}). The example below illustrates the technique for calculating the KLD between two Weibull densities with prescribed shape parameter.

\begin{example}
Consider the Weibull family with prescribed shape parameter $k$ that form an exponential family (including the family of exponential distributions for $k=1$ and the the Rayleigh distributions for $k=2$).
The density of a Weibull distribution with scale $\lambda>0$ and fixed shape parameter $k$ is
$$
p_\lambda(x)= \frac{k}{\lambda}\left(\frac{x}{\lambda}\right)^{k-1} e^{-(x / \lambda)^{k}},\quad x\in\calX=\bbR_{++}
$$
We have $\theta(u)=u^{-k}$ and $\theta'(u)=-ku^{-k-1}$.
We set $\omega=1$ so that $p_\lambda(\omega)=\frac{k}{\lambda^k}\exp(\lambda^{-k})$.
Let us program the formula of Eq.~\ref{eq:klefalpha} using the computer algebra system {\sc Maxima}:

\begin{verbatim}
/* KLD betwen Weibull distributions by calculating a limit */
declare( k , integer);
assume(lambda1>0);
assume(lambda2>0);
k:5;
omega:1;
t(u):=u**(-k);
tinv(u):=u**(-1/k);
tp(u):=k*u**(-k-1);
p(x,l):=(k/l)*((x/l)**(k-1))*exp(-(x/l)**k);
mean(a,b):= tinv(alpha*t(a)+(1-alpha)*t(b));
log(p(omega,l1)/p(omega,l2)) + (1.0/alpha)*log(p(omega,mean(l2,l1))/p(omega,l1));
limit (ratsimp(%), alpha, 0);
expand(%);
\end{verbatim}

\begin{figure}
\centering
\includegraphics[width=0.65\textwidth]{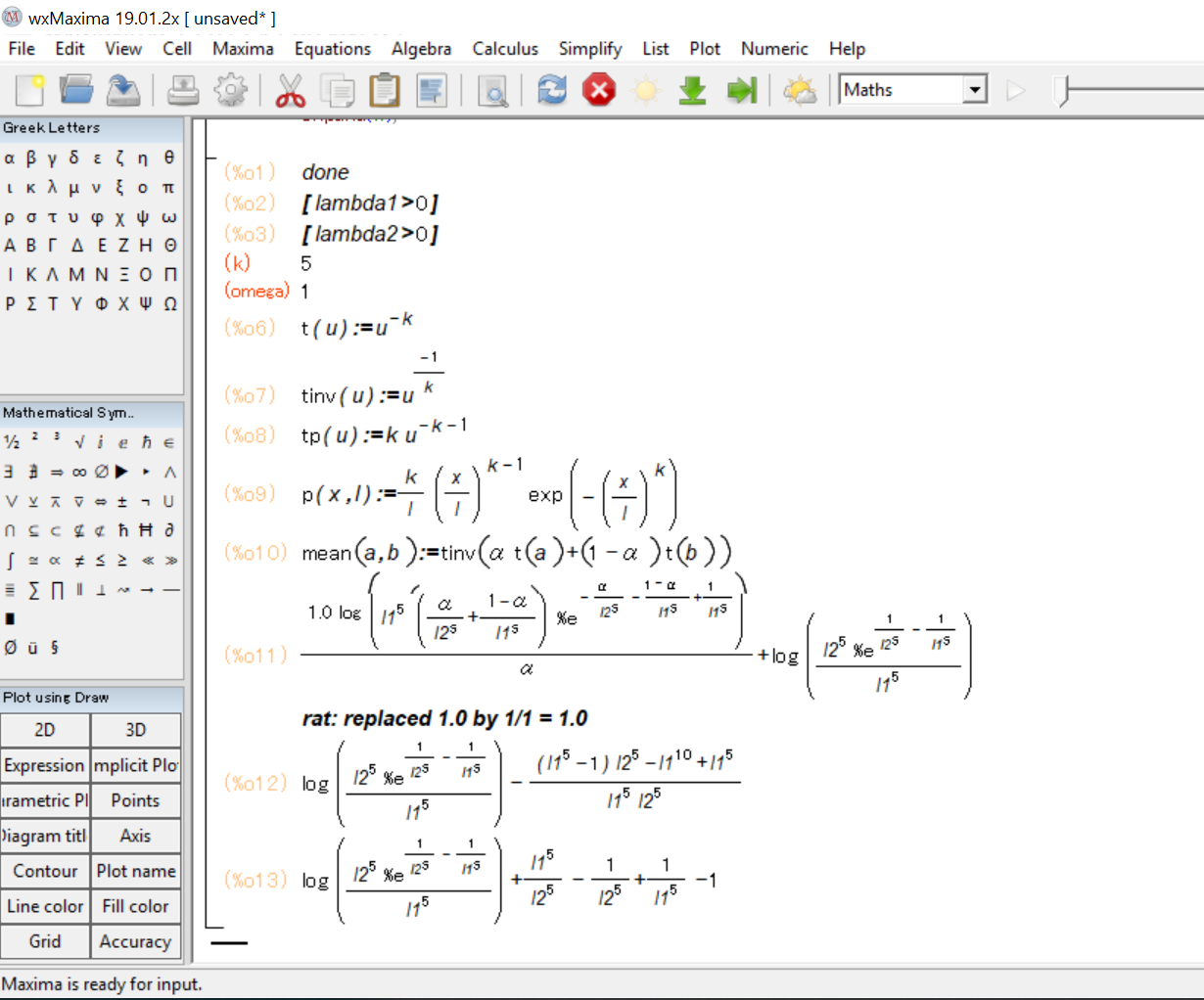}
\caption{Snapshot of the {\sc Maxima} GUI displaying the result of symbolic calculations of the KLD between two Weibull densities of prescribed parameter shape.\label{fig:weibullkllimit}}
\end{figure}

Figure~\ref{fig:weibullkllimit} displays a snapshot of the result which can be easily simplified manually as
\begin{equation}\label{eq:klweibullkfixed}
D_\KL[p_{\lambda_1}:p_{\lambda_2}] = k\log\frac{\lambda_2}{\lambda_1}+\left(\frac{\lambda_1}{\lambda_2}\right)^k-1.
\end{equation}

In general, the KLD between two Weibull densities with arbitrary shapes~\cite{KLWeibull-2013} is
\begin{equation}
D_\KL[p_{k_1,\lambda_1}:p_{k_2,\lambda_2}]= \log \frac{k_{1}}{\lambda_{1}^{k_{1}}}-\log \frac{k_{2}}{\lambda_{2}^{k_{2}}}+\left(k_{1}-k_{2}\right)\left(\log \lambda_{1}-
\frac{\gamma}{k_{1}}\right)+\left(\frac{\lambda_{1}}{\lambda_{2}}\right)^{k_{2}} \Gamma\left(\frac{k_{2}}{k_{1}}+1\right)-1,
\end{equation}
where $\gamma$ denotes the Euler-Mascheroni constant.
Thus when $k_1=k_2=k$, we recover Eq.~\ref{eq:klweibullkfixed} since $\Gamma(2)=1$.
(However, the family of Weibull densities with varying parameter shape is not an exponential family since the sufficient statistics depend on the shape parameters.)
\end{example}

In practice, we may program the formula of Eq.~\ref{eq:klefalpha} by defining:
\begin{equation}
D_{\KL,\alpha}[p_{\lambda_1}:p_{\lambda_2}] = \log\left(\frac{p(\omega;\lambda_1)}{p(\omega;\lambda_2)}\right) 
+   \frac{1}{\alpha} 
\log\left(\frac{p\left(\omega;M_{\theta,\alpha}(\lambda_2,\lambda_1)\right)}{p(\omega;\lambda_1)}\right),
\end{equation}
and approximate the KLD by $D_{\KL,\alpha}$ for a small value of $\alpha$ (say, $\alpha=10^{-3}$).
Thus we need only $\theta(\cdot)$ and $\theta^{-1}(\cdot)$ for defining $M_\theta(\cdot,\cdot)$, and the density $p(x;\theta)$.
This approximation also works for multivariate extensions of the quasi-arithmetic means.

Let us give some two examples using the first-order approximation of the univariate quasi-arithmetic mean: 

\begin{example}
Consider the family $\{p(x;\lambda)=\lambda\exp(-\lambda x)\}$ of exponential distributions with support $\calX=[0,\infty)$.
Set $\omega=0$, 
$p(\omega;\lambda)=\lambda$, $\theta=\lambda$, $\theta(u)=u$ and $\theta'(u)=1$.
We have
\begin{eqnarray}
D_\KL[p_{\lambda_1}:p_{\lambda_2}] &=& \log\left(\frac{\lambda_1}{\lambda_2}\right)
+ \lim_{\alpha\rightarrow 0} \frac{1}{\alpha} \log\frac{\lambda_1+\alpha(\lambda_2-\lambda_1)}{\lambda_1},\\
&=& \log\left(\frac{\lambda_2}{\lambda_1}\right)  +  
 \lim_{\alpha\rightarrow 0} \frac{1}{\alpha} \log \left(1+\alpha\left(\frac{\lambda_2}{\lambda_1}-1\right)\right),\\
&=&  \log\left(\frac{\lambda_1}{\lambda_2}\right)  + \frac{\lambda_2}{\lambda_1}-1. \label{eq:klexpdist}
\end{eqnarray}
\end{example}

\begin{example}
Consider the Poisson family with $\omega=0$,  $p_\lambda(\omega)=\exp(-\lambda)$, $\theta(u)=\log u$ ($M_\theta$ is the geometric mean) and $\theta'(u)=\frac{1}{u}$.
We get
\begin{eqnarray}
D_\KL[p_{\lambda_1}:p_{\lambda_2}] &=& \log \frac{p(\omega;\lambda_1)}{p(\omega;\lambda_2)} 
+ \lim_{\alpha\rightarrow 0} \frac{1}{\alpha} 
\log\frac{
p\left(\omega;\lambda_1+\alpha \frac{\theta(\lambda_2)-\theta(\lambda_1)
}{\theta'(\lambda_1)}\right)}{p(\omega;\lambda_2)},\\
&=& \lambda_2-\lambda_1 + \lim_{\alpha\rightarrow 0} \frac{1}{\alpha}  \log \exp(\alpha\lambda_1\log\frac{\lambda_1}{\lambda_2}),\\
&=& \lambda_2-\lambda_1 + \lambda_1\log\frac{\lambda_1}{\lambda_2}.
\end{eqnarray}
\end{example}

\subsection{Kullback-Leibler divergence formula relying on the differential entropy and moments}\label{sec:KLhmoment}
Consider the {\em Kullback-Leibler divergence}~\cite{KL-1951} (relative entropy) $D_\KL[p:q] = \int p(x)\log\frac{p(x)}{q(x)}\dmu(x)$ between two probability densities $p$ and $q$.
When the densities belong to the same exponential families, the KL divergences amounts to a {\em Legendre-Fenchel divergence} (the canonical expression of divergences using the dual coordinate systems in dually flat spaces of information geometry~\cite{IG-2016}):
\begin{equation}
D_\KL[p_{\lambda_1}:p_{\lambda_2}] = A_F(\theta_2:\theta_1),
\end{equation}
where the Legendre-Fenchel divergence is defined for a pair of strictly convex and differentiable conjugate generators $F(\theta)$ and $F^*(\eta)=\sup_{\theta} \theta^\top\eta-F(\theta)$ by
\begin{equation}\label{eq:bdip}
A_F(\theta:\eta') := F(\theta)+F^*(\eta')-\theta^\top \eta',
\end{equation}
with $\eta'=\nabla F(\theta)$.

Since $F$ is defined modulo some affine function, we can choose $F(\theta(\lambda))=-\log p(\omega;\theta(\lambda))$.
Furthermore, for exponential families, we have 
\begin{equation}
\eta(\lambda)=E_{p_{\lambda}}[t(x)],
\end{equation}
and
 the {\em Shannon entropy}~\cite{Shannon-1948} 
\begin{equation}
h(p)= -\int_{\calX} p(x)\log p(x) \dmu(x),
\end{equation}
admits the following expression~\cite{HEF-2010} when $p=p_\theta$ belongs to an exponential family $\calE$:
\begin{equation}
h(p_\lambda)=  -F^*(\eta(\lambda))-E_{p_\lambda}[k(x)].
\end{equation}

Thus if we already have at our disposal (1) the expectation of the sufficient statistics, and (2) the entropy, we can easily recover  the 
Kullback-Leibler divergence as follows:
\begin{equation}\label{eq:KLEFnocumulant}
D_\KL[p_{\lambda_1}:p_{\lambda_2}] = -\log p(\omega;\lambda_2) -h(p_{\lambda_1}) -E_{p_{\theta_1}}[k(x)] 
-\theta(\lambda_2)^\top E_{p_{\theta_1}}[t(x)].
\end{equation}

For densities $p_{\lambda_1}$ and $p_{\lambda_2}$ belonging to the same exponential family, the Jeffreys divergence is
\begin{equation}
D_J[p_{\lambda_1}:p_{\lambda_2}]=D_\KL[p:q]+D_\KL[q:p]=(\theta_2-\theta_1)^\top (\eta_2-\eta_1).
\end{equation}
It follows that
  we can write Jeffreys' divergence   using the following cumulant-free expression:
\begin{equation}
D_J[p_{\lambda_1}:p_{\lambda_2}]=\left(\theta(\lambda_2)-\theta(\lambda_1)\right)^\top (E_{p_{\lambda_2}}[t(x)] - E_{p_{\lambda_1}}[t(x)]).
\end{equation}

Note that a strictly monotone operator $O$ defines a symmetric dissimilarity: 
\begin{equation}
D_O(\theta_1,\theta_2):=(\theta_1-\theta_2)^\top (O(\theta_2)-O(\theta_1))\geq 0,
\end{equation}
with equality if and only if $\theta_1=\theta_2$.
Since $\nabla F$ is a strictly monotone operator and $E_{p_{\lambda}}[k(x)]=\nabla F(\theta)$, we may reinterpret 
the Jeffreys' divergence as a symmetric dissimilarity induced by a strictly monotone operator.

Let us report now some illustrating examples.
We start with an example illustrating the use of a separable  multivariate quasi-arithmetic mean.

\begin{example}[continue Example~\ref{ex:DirichletBhat}]
Consider the Dirichlet exponential family.
The differential entropy of a Dirichlet density $p_\alpha$ is 
$$
h(p_{\alpha})=\log \mathrm{B}(\alpha)+\left(\sum_i\alpha_{i}-d\right) \psi\left(\sum_i\alpha_{i}\right)-\sum_{j=1}^{d}\left(\alpha_{j}-1\right) \psi\left(\alpha_{j}\right),
$$
where $\psi(\cdot)$ denotes the  digamma function.
We have $E[t(x)]=E[(\log x_1,\ldots,\log x_d)]=(\psi(\alpha_1)-\psi(\alphasum), \ldots, \psi(\alpha_1)-\psi(\alphasum))$ 
where $\alphasum=\sum_{i=1}^d \alpha_i$.
It follows that the Kullback-Leibler between $p_{\alpha}$ and $p_{\alpha'}$ is:
$$
D_\KL[p_{\alpha}:p_{\alpha}'] = \log \Gamma\left(\alpha_{\sum}\right)-\sum_{i=1}^{d} \log \Gamma\left(\alpha_{i}\right)-\log \Gamma\left(\alpha_{\sum}'\right)+\sum_{i=1}^{d} \log \Gamma\left(\alpha_{i}'\right)+\sum_{i=1}^{d}\left(\alpha_{i}-\alpha_{i}'\right)
\left(\psi\left(\alpha_{i}\right)-\psi\left(\alpha_{\sum}\right)\right).
$$
\end{example}

Next, we report an example illustrating a non-separable multivariate quasi-arithmetic mean.

\begin{example}[continue Example~\ref{ex:GaussianBhat}]
Consider the $d$-dimensional multivariate Gaussian family.
Since $k(x)=0$ (so that $E[k(x)]=0$), 
$\eta(\lambda)=\nabla_\theta F(\theta(\lambda))=(\mu,\mu\mu^\top+\Sigma)$
 (since $E[x]=\mu$, $E[xx^\top]=\mu\mu^\top+\Sigma$), and the usual differential entropy is known  as
$h(p_{\mu,\Sigma})=\frac{1}{2}\log |2\pi e\Sigma|= \frac{d}{2}(1+\log 2\pi)+\frac{1}{2}\log |\Sigma|$ since $|aM|=|a|^d |M|$.
we recover the Kullback-Leibler divergence as
\begin{equation}\label{eq:klmvn}
D_\KL[p_{\mu_1,\Sigma_1}:p_{\mu_2,\Sigma_2}] =
\frac{1}{2} \left( \tr(\Sigma_2^{-1}\Sigma_1) + \Delta_\mu^\top \Sigma_2^{-1}\Delta_\mu +  \log \frac{|\Sigma_2|}{|\Sigma_1|} -d \right), 
\end{equation}
and the  Jeffreys divergence as 
\begin{equation}
D_J[p_{\mu_1,\Sigma_1},p_{\mu_2,\Sigma_2}] = 
     \Delta_\mu^\top \left(\frac{\Sigma_1^{-1}+\Sigma_2^{-1}}{2}\right) \Delta_\mu 
		+ \frac{1}{2} \tr\left(\Sigma_2^{-1}\Sigma_1+\Sigma_1^{-1}\Sigma_2\right)
		- d.
\end{equation}
\end{example}

\subsection{The Kullback-Leibler divergence expressed as a log density ratio}\label{sec:KLdr}

Let us express the Bregman divergence with the {\em equivalent} generator $F_l(\theta) :=  -\log(p_\theta(\omega))$ for any 
prescribed $\omega\in\calX$ instead of the cumulant function of Eq.~\ref{eq:cumulantF}.
We get
\begin{eqnarray}\label{eq:klomega}
D_\KL[p_{\lambda_1}:p_{\lambda_2}] &=& B_{F_l}(\theta(\lambda_2):\theta(\lambda_1)),\\
&=& \log\left( \frac{p_{\lambda_1}(\omega)}{p_{\lambda_2}(\omega)} \right)
 -(\theta(\lambda_2)-\theta(\lambda_1))^\top \nabla_\theta F_l(\theta(\lambda_1)).
\end{eqnarray}
Let us remark that 
\begin{equation}\label{eq:Fl}
F_l(\theta)=-\theta^\top t(\omega)+F(\theta)-k(\omega).
\end{equation}
We have
\begin{eqnarray}
\nabla  F_l(\theta(\lambda_1)) &=& - \frac{\nabla_\theta p_{\lambda_1}(\omega)}{p_{\lambda_1}(\omega)},\\
&=& - \frac{(t(\omega)-\nabla F(\theta(\lambda_1)))p_{\lambda_1}(\omega)}{p_{\lambda_1}(\omega)},\\
&=&  - (t(\omega)-\nabla F(\theta(\lambda_1))),
\end{eqnarray}
where $F(\theta)$ is the cumulant function of Eq.~\ref{eq:cumulantF}.
Alternatively, we can also directly find that $\nabla  F_l(\theta(\lambda_1))=- (t(\omega)-\nabla F(\theta(\lambda_1)))$ from Eq.~\ref{eq:Fl}.

It follows that:
\begin{equation}\label{eq:klomegaf}
D_\KL[p_{\lambda_1}:p_{\lambda_2}] = 
\log\left( \frac{p_{\lambda_1}(\omega)}{p_{\lambda_2}(\omega)} \right)
+ (\theta(\lambda_2)-\theta(\lambda_1))^\top (t(\omega)-\nabla F(\theta(\lambda_1))),\quad
\forall \omega\in\calX.
\end{equation}

Thus when $t(\omega)-\nabla F(\theta(\lambda_1))$ is Euclidean orthogonal to $\theta(\lambda_2)-\theta(\lambda_1)$  (i.e., $ (\theta(\lambda_2)-\theta(\lambda_1))^\top (t(\omega)-\nabla F(\theta(\lambda_1)))=0$), the Bregman divergence (and the corresponding KLD on swapped parameter order of the densities) is expressed as a log-density ratio quantity.
Let 
\begin{equation}
\calX_\perp^\calE(\lambda_1:\lambda_2) :=\left\{\omega\in\calX\ :\ (\theta(\lambda_2)-\theta(\lambda_1))^\top (t(\omega)-\nabla F(\theta(\lambda_1))) =0\right\}.
\end{equation}
Then $D_\KL[p_{\lambda_1}:p_{\lambda_2}] =\log\left( \frac{p_{\lambda_1}(\omega)}{p_{\lambda_2}(\omega)} \right)$ 
for all $\omega\in \calX_\perp^\calE(\lambda_1:\lambda_2)$.

\begin{lemma}\label{thm:klratio}
The Kullback-Leibler divergence between two  densities $p_{\lambda_1}$ and 
$p_{\lambda_2}$ belonging to a same exponential family $\calE$ is expressed as a log density ratio,
 $D_\KL[p_{\lambda_1}:p_{\lambda_2}] =\log\left( \frac{p_{\lambda_1}(\omega)}{p_{\lambda_2}(\omega)} \right)$, whenever
$\omega\in \calX_\perp^\calE(\lambda_1:\lambda_2)$.
\end{lemma} 

Notice that a sufficient condition is to choose $\omega$ such that $t(\omega)=\eta_1=\nabla F(\theta_1)$.

Thus if we carefully choose $\omega\in \calX_\perp^\calE(\lambda_1:\lambda_2)$ according to the source parameters, we may express the Kullback-Leibler divergence as a simple log density ratio without requiring the formula for the differential entropy nor the moment.

\begin{example}
Consider the exponential family $\calE=\{p_\lambda(x)=\lambda\exp(-\lambda x),\ \lambda>0\}$ of exponential distributions with $\nabla F(\theta)=\frac{1}{\theta}$ for $\theta=\lambda$.
We have $(\theta(\lambda_2)-\theta(\lambda_1))^\top (t(\omega)-\nabla F(\theta(\lambda_1)))=0$ that amounts to 
$(\lambda_2-\lambda_1)(\omega-\frac{1}{\lambda_1})=0$, i.e., $\omega=\frac{1}{\lambda_1}$ 
(and $\calX_\perp^\calE(\lambda_1:\lambda_2)=\left\{\frac{1}{\lambda_1}\right\}$).
In that case, we have 
\begin{eqnarray}
D_\KL[p_{\lambda_1}:p_{\lambda_2}] &=& \log\left( \frac{p_{\lambda_1}\left(\frac{1}{\lambda_1}\right)}{p_{\lambda_2}\left(\frac{1}{\lambda_1}\right)} \right),\\
&=& \log \left(
\frac{ \lambda_1 \exp\left(-\lambda_1 \frac{1}{\lambda_1}\right) }{\lambda_2 \exp\left(-\lambda_2 \frac{1}{\lambda_1}\right)}
\right),\\
&=& \log\left( \frac{\lambda_1}{\lambda_2}\right) + \frac{\lambda_2}{\lambda_1}-1.
\end{eqnarray}
This formula matches the expression of Eq.~\ref{eq:klexpdist}.
\end{example}

Similarly, we may rewrite the Bregman divergence $B_F(\theta_1:\theta_2)$ as
\begin{equation}
B_F(\theta_1:\theta_2)=B_{F_a}(\theta_1:\theta_2)=F_a(\theta_1)-F_a(\theta_2),
\end{equation}
for $F_a(\theta)=F(\theta)-a\theta$ (and $\nabla F_a(\theta)=\nabla F(\theta)-a$) for $a=-\nabla F_a(\theta_2)$.

\begin{example}
Consider the exponential family of $d$-dimensional Gaussians with fixed covariance matrix $\Sigma$:
\begin{equation}
\calE=\left\{p_\lambda(x)=  \frac{1}{(2\pi)^{\frac{d}{2}}\sqrt{|\Sigma|}}  \exp\left(-\frac{1}{2} (x-\lambda)^\top \Sigma^{-1} (x-\lambda)\right), \ :\ \lambda\in\bbR^d \right\}.
\end{equation}
It is an exponential family of order $D=d$ with sufficient statistic $t(x)=x$.
Let us choose $\omega$ such that $(\theta(\lambda_2)-\theta(\lambda_1))^\top (t(\omega)-\nabla F(\theta(\lambda_1))) =0$.
For example, we choose $t(\omega)=\nabla F(\theta_1)=\mu_1$.
It follows that we have
\begin{eqnarray}
D_\KL[p_{\lambda_1}:p_{\lambda_2}] &=& \log\left( \frac{p_{\lambda_1}\left(\frac{1}{\lambda_1}\right)}{p_{\lambda_2}\left(\frac{1}{\lambda_1}\right)} \right),\\
&=& \log\left(
\frac{(2\pi)^{\frac{d}{2}}\sqrt{|\Sigma|}}{(2\pi)^{\frac{d}{2}}\sqrt{|\Sigma|}}
\exp\left(
\frac{1}{2} (\mu_1-\mu_2)^\top \Sigma^{-1} (\mu_1-\mu_2)
\right)
\right),\\
&=& \frac{1}{2} (\mu_1-\mu_2)^\top \Sigma^{-1} (\mu_1-\mu_2).
\end{eqnarray}
This is half of the squared Mahalanobis distance obtained for the precision matrix $\Sigma^{-1}$.
We recover the Kullback-Leibler divergence between multivariate Gaussians (Eq.~\ref{eq:klmvn}) when $\Sigma_1=\Sigma_2=\Sigma$.
\end{example}

\begin{example}
Consider the exponential family of Rayleigh distributions:
\begin{equation}
\calE=\left\{ p_\lambda(x)= \frac{x}{\lambda^{2}} \exp \left(-\frac{x^{2}}{2 \lambda^{2}}\right) \right\},
\end{equation}
for $\calX=[0,\infty)$.
Let us choose $\omega$ such that $(\theta(\lambda_2)-\theta(\lambda_1))^\top (t(\omega)-\nabla F(\theta(\lambda_1))) =0$ 
with $\theta=-\frac{1}{2\lambda^2}$, $t(x)=x^2$ and $\nabla F(\theta(\lambda))=2\lambda^2$.
We choose $\omega^2=2\lambda^2_1$ (i.e., $\omega=\lambda_1\sqrt{2}$).
It follows that we have
\begin{eqnarray}
D_\KL[p_{\lambda_1}:p_{\lambda_2}] &=& \log\left( \frac{p_{\lambda_1}\left(\frac{1}{\lambda_1}\right)}{p_{\lambda_2}\left(\frac{1}{\lambda_1}\right)} \right),\\
&=& \log\left(
\frac{\sqrt{2}}{\lambda_1}\frac{\lambda_2^2}{\lambda_1}\exp\left(
-\frac{2\lambda_1^2}{2\lambda_1^2}+\frac{2\lambda_1^2}{2\lambda_2^2}
\right)
 \right),\\
&=& 2\log\left(\frac{\lambda_2}{\lambda_1}\right) + \frac{\lambda_1^2}{\lambda_2^2} -1.]
\end{eqnarray}  
\end{example}

This example shows the limitation of the method which we shall now overcome.

\begin{example}\label{ex:klgaussian}
Consider the exponential family of univariate Gaussian distributions:
\begin{equation}
\calE=\left\{ p_\lambda(x)= \frac{1}{\sqrt{2 \pi \lambda_2}} \exp \left(-\frac{(x-\lambda_1)^{2}}{2 \lambda_2}\right)  \right\},
\end{equation}
for $\calX=(-\infty,\infty)$.
Let us choose $\omega$ such that $(\theta(\lambda_2)-\theta(\lambda_1))^\top (t(\omega)-\nabla F(\theta(\lambda_1))) =0$.
Here, $\nabla F(\theta(\lambda_1)))=(\mu_1,\mu_1^2+\sigma_1^2)$ and $t(\omega)=(x,x^2)$.
Thus we have $t(\omega)\not=\nabla F(\theta(\lambda_1))$ for any $\omega\in\calX$.
\end{example}

Let $H=\{\nabla F(\theta)\ :\ \theta\in\Theta\}$ denote the dual moment parameter space.
When the exponential family is regular, $\Theta$ is an open convex set, and so is $H$ ($F$ is of Legendre-type).
The problem we faced with the last example is that $\omega\in\partial H$.
However since Eq.~\ref{eq:klomegaf} holds for any $\omega\in\Omega$, let us choose $s$ values for $\omega$ (i.e., $\omega_1,\ldots,\omega_s$), and average Eq.~\ref{eq:klomegaf} for these $s$ values.
We have
\begin{eqnarray}\label{eq:klomegaff}
D_\KL[p_{\lambda_1}:p_{\lambda_2}] &=& \frac{1}{s} \sum_{i=1}^s  
\log\left( \frac{p_{\lambda_1}(\omega_i)}{p_{\lambda_2}(\omega_i)} \right)
+ (\theta(\lambda_2)-\theta(\lambda_1))^\top \left(\frac{1}{s}\sum_{i=1}^s t(\omega_i)-\nabla F(\theta(\lambda_1))\right).
\end{eqnarray}

We can now choose the $\omega_i$'s such that $\frac{1}{s}\sum_{i=1}^s t(\omega_i)\in H$ for $s>1$.
We need to choose $s$ so that the system of equations $\frac{1}{s}t(\omega_i)=\nabla F(\theta)=E[t(x)]$ is solvable.

\begin{example}
Let us continue Example~\ref{ex:klgaussian}.
Consider $s=2$. We need to find $\omega_1$ and $\omega_2$ such that we can solve the following system of equations:
\begin{equation}
\left\{
\begin{array}{lcl}
\mu_1&=&\frac{\omega_1+\omega_2}{2},\\
\mu_1^2+\sigma_1^2&=&\frac{\omega_1^2+\omega_2^2}{2}.
\end{array}
\right.
\end{equation}

The solution of this system of equations is $\omega_1=\mu_1-\sigma_1$ and $\omega_2=\mu_1+\sigma_1$.
Thus it follows that we have:
\begin{equation}
D_\KL[p_{\mu_1,\sigma_1}:p_{\mu_2,\sigma_2}] = \frac{1}{2}
\left(    
 \log\left( \frac{p_{\mu_1,\sigma_1}(\mu_1-\sigma_1)}{p_{\mu_2,\sigma_2}(\mu_1-\sigma_1)} \right) 
+\log\left( \frac{p_{\mu_1,\sigma_1}(\mu_1+\sigma_1)}{p_{\mu_2,\sigma_2}(\mu_1+\sigma_1)} \right)
\right).
\end{equation}
\end{example}

We conclude with the following theorem extending Lemma~\ref{thm:klratio}:

\begin{theorem}\label{thm:klratiofull}
The Kullback-Leibler divergence between two  densities $p_{\lambda_1}$ and 
$p_{\lambda_2}$ belonging to a full regular exponential family $\calE$ of order $D$   can be expressed as the averaged sum of logarithms of density ratios:
 $$
D_\KL[p_{\lambda_1}:p_{\lambda_2}] = \frac{1}{s}\sum_{i=1}^s \log\left( \frac{p_{\lambda_1}(\omega_i)}{p_{\lambda_2}(\omega_i)} \right),
$$
 where
$\omega_1,\ldots, \omega_s$ are $s\leq D+1$ distinct samples of $\calX$ chosen such that 
$\frac{1}{s}\sum_{i=1}^s t(\omega_i)=E_{p_{\lambda_1}}[t(x)]$.
\end{theorem}

The bound $s\leq D+1$ follows from  Carath\'eodory's theorem~\cite{Caratheodory-1907}.

Notice that the Monte Carlo stochastic approximation of the Kullback-Leibler divergence:
\begin{equation}
\widetilde{D}_\KL^{(m)}[p_{\lambda_1}:p_{\lambda_2}] = \frac{1}{m} \sum_{i=1}^m \log\frac{p_{\lambda_1}(x_i)}{p_{\lambda_2}(x_i)},
\end{equation}
for $m$ independent and identically distributed (iid) variates $x_1,\ldots, x_m$  from $p_{\lambda_1}$ is also an average sum of log density ratios which yields a consistent estimator, i.e., 
\begin{equation}
\lim_{m\rightarrow\infty} \widetilde{D}_\KL^{(m)}[p_{\lambda_1}:p_{\lambda_2}]=D_\KL[p_{\lambda_1}:p_{\lambda_2}].
\end{equation}

In practice, to avoid potential negative values of $\widetilde{D}_\KL^{(m)}[p_{\lambda_1}:p_{\lambda_2}]$, we estimate the extended Kullback-Leibler divergence:
\begin{equation}
\widetilde{D}_{\KL_+}^{(m)}[p_{\lambda_1}:p_{\lambda_2}] = \frac{1}{m} \sum_{i=1}^m \left(\log\frac{p_{\lambda_1}(x_i)}{p_{\lambda_2}(x_i)}+p_{\lambda_2}(x_i)-p_{\lambda_1}(x_i)\right)\geq 0,
\end{equation}
with $\lim_{m\rightarrow\infty} \widetilde{D}_{\KL_+}^{(m)}[p_{\lambda_1}:p_{\lambda_2}]=D_\KL[p_{\lambda_1}:p_{\lambda_2}]$.

\begin{example}
We continue Example~\ref{ex:GaussianBhat} of the $d$-dimensional multivariate Gaussian family.
First, let us  consider the subfamily of zero-centered Gaussian densities. 
The sufficient statistic $t(x)$ is a $d\times d$ matrix: $t(x)=xx^\top$.
We find the $d$ column vectors $\omega_i$'s from the singular value decomposition of $\Sigma_1$:
\begin{equation}
\Sigma_1 = \sum_{i=1}^d \lambda_i e_ie_i^\top,
\end{equation}
where the $\lambda_i$'s are the eigenvalues and the $e_i$'s the corresponding eigenvectors.
Let $\omega_i=\sqrt{d\lambda_i}e_i$ for $i\in\{1,\ldots, d\}$.
We have
\begin{equation}
\frac{1}{d}\sum_{i=1}^d t(\omega_i)=\Sigma_1=E_{p_{\Sigma_1}}[xx^\top].
\end{equation}

It follows that we can express the KLD between two zero-centered Gaussians $p_{\Sigma_1}$ and $p_{\Sigma_2}$
 as the following weighted sum of log density ratios:
\begin{equation}
D_\KL[p_{\Sigma_1}:p_{\Sigma_2}] = 
\frac{1}{d} \sum_{i=1}^d \log\left( \frac{p_{\Sigma_1}(\sqrt{d\lambda_i}e_i)}{p_{\Sigma_2}(\sqrt{d\lambda_i}e_i)} \right),
\end{equation}
where  the $\lambda_i$'s are the eigenvalues of $\Sigma_1$ with the $e_i$'s the corresponding eigenvectors.
Notice that the order of the family is $D=\frac{d(d+1)}{2}$ but we used $s=d\leq D$ vectors $\omega_1,\ldots,\omega_d$.

The closed-form formula for the Kullback-Leibler divergence between two zero-centered Gaussians 
$p_{\Sigma_1}$ and $p_{\Sigma_2}$ is 
\begin{equation}
D_\KL[p_{\Sigma_1}:p_{\Sigma_2}] = \frac{1}{2}\left(
\tr\left(\Sigma_2\Sigma_1^{-1}\right)+\log\left(\frac{|\Sigma_2|}{|\Sigma_1|}\right)-d
\right).
\end{equation}

Now consider the full family of multivariate normal densities.
We shall use  $2d$ vectors $\omega_i$ as follows:
\begin{eqnarray}
\left\{
\begin{array}{ll}
\omega_i=\mu_1-\sqrt{d\lambda_i}e_i, & i\in\{1,\ldots, d\}\\
\omega_i=\mu_1+\sqrt{d\lambda_i}e_i, & i\in\{d+1,\ldots, 2d\}
\end{array}
\right.
\end{eqnarray}
where  the $\lambda_i$'s are the eigenvalues of $\Sigma_1$ with the $e_i$'s the corresponding eigenvectors.
It can be checked that $\frac{1}{2d}\sum_{i=1}^{2d} \omega_i= E_{p_{\lambda_1}}[x]=\mu_1$
 and $\frac{1}{2d}\sum_{i=1}^{2d} \omega_i\omega_i^\top=\mu_1\mu_1^\top+\Sigma_1$.
These points are called the ``sigma points'' in the unscented transform~\cite{julier2000new,goldberger2008simplifying}.
Moreover, we have $\sqrt{\lambda_i}e_i=[\sqrt{\Sigma_1}]_{\cdot,i}$, the $i$-th column of the square root of the covariance matrix of $\Sigma_1$.
Thus it follows that
\begin{eqnarray*}
D_\KL[p_{\mu_1,\Sigma_1}:p_{\mu_2,\Sigma_2}] &=& 
\frac{1}{2d} \sum_{i=1}^{d}
\left( \log\left( 
\frac{p_{\mu_1,\Sigma_1}\left(\mu_1-\sqrt{d\lambda_i}e_i\right)}{p_{\mu_2,\Sigma_2}\left(\mu_1-\sqrt{d\lambda_i}e_i\right) }
\right)
+
\log\left( 
\frac{p_{\mu_1,\Sigma_1}\left(\mu_1+\sqrt{d\lambda_i}e_i\right)}{p_{\mu_2,\Sigma_2}\left(\mu_1+\sqrt{d\lambda_i}e_i\right)} 
  \right)
	\right),\\
	&=&
	\frac{1}{2d} \sum_{i=1}^{d}
\left( \log\left( 
\frac{p_{\mu_1,\Sigma_1}\left(\mu_1-[\sqrt{d\Sigma_1}]_{\cdot,i}\right)}{p_{\mu_2,\Sigma_2}\left(\mu_1-[\sqrt{d\Sigma_1}]_{\cdot,i}\right) }
\right)
+
\log\left( 
\frac{p_{\mu_1,\Sigma_1}\left(\mu_1+[\sqrt{d\Sigma_1}]_{\cdot,i}\right)}{p_{\mu_2,\Sigma_2}\left(\mu_1+[\sqrt{d\Sigma_1}]_{\cdot,i}\right)} 
  \right)
	\right),
\end{eqnarray*}
where $[\sqrt{d\Sigma_1}]_{\cdot,i}=\sqrt{\lambda_i}e_i$ denotes the vector extracted from the $i$-th column of the square root matrix of $d\Sigma_1$.

This formula matches the ordinary formula for the Kullback-Leibler divergence between the two Gaussian densities 
$p_{\mu_1,\Sigma_1}$ and $p_{\mu_2,\Sigma_2}$:
\begin{equation}
D_\KL[p_{\mu_1,\Sigma_1}:p_{\mu_2,\Sigma_2}] = 
\frac{1}{2}\left(  (\mu_2-\mu_1)^\top\Sigma_2^{-1}(\mu_2-\mu_1)+
\tr\left(\Sigma_2\Sigma_1^{-1}\right)+\log\left(\frac{|\Sigma_2|}{|\Sigma_1|}\right)-d
\right).
\end{equation}			
\end{example}

\begin{example}
We continue the Example~\ref{ex:Gamma} of the exponential family of gamma distributions which has order $D=2$.
The sufficient statistic vector is $t(x)=(x,\log x)$, and we have $E[x]=\frac{\alpha}{\beta}$ and $E[\log x]=\psi(\alpha)-\log\beta$, where
 $\psi(\cdot)$ is the digamma function.
We want to express $D_\KL[p_{\alpha_1,\beta_1}:p_{\alpha_2,\beta_2}]$ as an average sum of log ratio of densities.
To find the values of $\omega_1$ and $\omega_2$, we need to solve the following system of equations:
\begin{equation}
\left\{
\begin{array}{lll}
\frac{\omega_1+\omega_2}{2}&=&a,\\
\frac{\log\omega_1+\log\omega_2}{2}&=& b
\end{array}
\right.,
\end{equation}
with $a=\frac{\alpha_1}{\beta_1}$ and $b=\psi(\alpha_1)-\log\beta_1$.
We find the following two solutions:
\begin{equation}
\omega_1=a-\sqrt{a^2-\exp(2b)},\quad \omega_2=a+\sqrt{a^2-\exp(2b)}.
\end{equation}
We have $a^2-\exp(2b)\geq 0$ since $E[x]^2\leq \exp(2E[\log x])$.

It follows that
\begin{eqnarray}
D_\KL[p_{\alpha_1,\beta_1}:p_{\alpha_2,\beta_2}] &=& \frac{1}{2}\left(\log\frac{p_{\alpha_1,\beta_1}(\omega_1)}{p_{\alpha_2,\beta_2}(\omega_1)}  + \log\frac{p_{\alpha_1,\beta_1}(\omega_2)}{p_{\alpha_2,\beta_2}(\omega_2)}  \right),
\end{eqnarray}
This expression of the KLD is the same as the ordinary expression of the KLD:
\begin{equation}
D_\KL[p_{\alpha_1,\beta_1}:p_{\alpha_2,\beta_2}]=(\alpha_1-\alpha_2)\psi(\alpha_1)-\log\Gamma(\alpha_1)+\log\Gamma(\alpha_2)+\alpha_2\log\frac{\beta_1}{\beta_2}+\alpha_1\frac{\beta_2-\beta_1}{\beta_1}.
\end{equation}
\end{example}

\begin{example}
Consider the exponential family of Beta distributions:
\begin{equation}
\calE=\left\{
p_{\alpha,\beta}(x) = \frac{x^{\alpha-1}(1-x)^{\beta-1}}{B(\alpha,\beta)}\ : \alpha>0,\beta>0, x\in\calX=(0,1)
\right\},
\end{equation}
where
\begin{equation}
B(\alpha,\beta)=\frac{\Gamma(\alpha)\Gamma(\beta)}{\Gamma(\alpha+\beta)}
\end{equation}
The sufficient statistic vector is $t(x)=(\log x,\log(1-x))$.
We have $E_{p_{\alpha_1,\beta_1}}[\log(x)]=\psi(\alpha_1)-\psi(\alpha_1+\beta_1)=:A$
and $E_{p_{\alpha_1,\beta_1}}[\log(1-x)]=\psi(\beta_1)-\psi(\alpha_1+\beta_1)=:B$.
We need to solve the following system of equations for $s=2$:
\begin{equation}
\left\{
\begin{array}{lcl}
\frac{\log\omega_1+\log\omega_2}{2} &=& A,\\
\frac{\log(1-\omega_1)+\log(1-\omega_2)}{2} &=& B
\end{array}
\right.
\end{equation}
This system is equivalent to
\begin{equation}
\left\{
\begin{array}{lcl}
\omega_1\omega &=& \exp(2A):=a,\\
(1-\omega_1)(1-\omega_2) &=& \exp(2B):=b
\end{array}
\right.
\end{equation}
We find
\begin{eqnarray}
\omega_1 &=& \frac{-b+a+1-\sqrt{\Delta}}{2},\\
\omega_2 &=& \frac{-b+a+1+\sqrt{\Delta}}{2},
\end{eqnarray}
where
\begin{equation}
\Delta=b^2-2(a+1)b+a^2-2a+1.
\end{equation}
It follows that we have
\begin{equation}
D_\KL[p_{\alpha_1,\beta_1}:p_{\alpha_2,\beta_2}] 
=
\frac{1}{2}\left(
\log\frac{p_{\alpha_1,\beta_1}(\omega_1)}{p_{\alpha_2,\beta_2}(\omega_1)} +
\log\frac{p_{\alpha_1,\beta_1}(\omega_2)}{p_{\alpha_2,\beta_2}(\omega_2)} 
\right).
\end{equation}

This formula is equivalent to the ordinary formula for the KLD between two beta densities $p_{\alpha_1,\beta_1}$ and $p_{\alpha_2,\beta_2}$:
\begin{eqnarray}
\lefteqn{D_\KL[p_{\alpha_1,\beta_1}:p_{\alpha_2,\beta_2}]=}\nonumber\\
&&\log \left(\frac{\mathrm{B}\left(\alpha_2, \beta_2\right)}{\mathrm{B}(\alpha_1, \beta_1)}\right)+\left(\alpha_1-\alpha_2\right) \psi(\alpha_1)+\left(\beta_1-\beta_2\right) \psi(\beta_1)+\left(\alpha_2-\alpha_1+\beta_2-\beta_1\right) \psi(\alpha_1+\beta_1)\nonumber.\\
\end{eqnarray}

\end{example}

Notice that the $\omega_i$'s are chosen according to $\lambda_1$.
Thus we may express the Voronoi bisector: 
\begin{equation}
\mathrm{Bi}(p_{\lambda_1},p_{\lambda_2}) := \left\{\lambda\ :\ D_\KL[p_{\lambda}:p_{\lambda_1}]=D_\KL[p_{\lambda}:p_{\lambda_2}] \right\}
\end{equation}
as follows:
\begin{equation}
\mathrm{Bi}(p_{\lambda_1},p_{\lambda_2})=\left\{\lambda\ :\ \frac{1}{s}\sum_{i=1}^s \log\left( \frac{p_{\lambda}(\omega_i)}{p_{\lambda_1}(\omega_i)}\right)  
= \frac{1}{s}\sum_{i=1}^s \log\left( \frac{p_{\lambda}(\omega_i)}{p_{\lambda_2}(\omega_i)}\right)
\right\}.
\end{equation}

In particular, when $s=1$, the Voronoi bisector is expressed as:
\begin{equation}
\mathrm{Bi}(p_{\lambda_1},p_{\lambda_2})=\left\{\lambda\ :\  p_{\lambda_1}(\omega(\lambda))=p_{\lambda_2}(\omega(\lambda))
\right\},
\end{equation}
where $\omega(\lambda)$ emphasizes on the fact that $\omega$ is a function of $\lambda$.
The statistical Voronoi diagrams of a finite set of densities belonging to an exponential family has been studied as equivalent Bregman Voronoi diagrams in~\cite{BVD-2010}.

\subsection{The Jensen-Shannon divergence}

The Jensen-Shannon divergence~\cite{JS-2019} (JSD) is another symmetrization of the Kullback-Leibler divergence which can be given many information-theoretic interpretations~\cite{JSD-2020} and which is further guaranteed to be always 
bounded by $\log 2$ (KLD and JD are unbounded):
\begin{eqnarray}
D_\JS[p,q] &:=& \frac{1}{2} \left( D_\KL\left[p:\frac{p+q}{2}\right] +  D_\KL\left[q:\frac{p+q}{2}\right] \right),\\
&=& h\left(\frac{p+q}{2}\right)-\frac{h(p)+h(q)}{2}.\label{eq:jsdh}
\end{eqnarray}

Usually, the JSD does not provably admit a closed-form formula~\cite{MixtureDiv-2016}.
However, in the particular case when the mixture $\frac{p_{\theta_1}+p_{\theta_2}}{2}$ belongs to the same parametric family of densities, we can calculate the Jensen-Shannon divergence using the entropy function as shown in Eq.~\ref{eq:jsdh}.
For example, consider a mixture family in information geometry~\cite{IG-2016}.
That is, a statistical mixture with $k$ prescribed  components $p_1(x), \ldots, p_k(x)$  which are linearly independent (so that all mixtures of the family are identifiable by their corresponding parameters).
Let $m_\lambda(x)=\sum_{i=1}^k w_ip_i(x)$.
In that particular case (e.g., mixture family with $k$ prescribed Gaussians components), we get
\begin{equation}
\frac{m_{\lambda_1}+m_{\lambda_2}}{2}=m_{\frac{\lambda_1+\lambda_2}{2}}.
\end{equation}
Thus the JSD for a mixture family can be expressed using the entropy as:
\begin{equation}
D_\JS[m_{\lambda_1},m_{\lambda_1}]= h\left(m_{\frac{\lambda_1+\lambda_2}{2}}\right) -\frac{h(m_{\lambda_1})+h(m_{\lambda_2})}{2}.
\end{equation}

Although we do not have closed-form formula for the entropy of a mixture (except in few cases, e.g., when the support of the distributions are pairwise disjoint~\cite{JSD-2020}), but we can use any approximation method for calculating the entropy of a mixture to approximate or bound~\cite{MixtureDiv-2016} the Jensen-Shannon divergence $D_\JS$.

\section{Conclusion}\label{sec:concl}

We have described several methods to easily recover closed-form formula for some common (dis)similarities between densities belonging to a same exponential family $\{p(x;\lambda)\}_{\lambda\in\Lambda}$ which express themselves using the cumulant function $F$ of the exponential family (e.g., the Kullback-Leibler divergence amounts to a reverse Bregman divergence and the Bhattacharyya distance amounts to a Jensen divergence).
Our trick consists in observing that the generators $F$ of the Bregman or Jensen divergences are defined modulo an affine term, so that we may choose $F(\theta(\lambda))=-\log p(\omega,\lambda)$ for {\em any} $\omega$ falling inside the support $\calX$.  
It follows that the Bhattacharyya coefficient can be calculated with the following {\em cumulant-free} expression:
\begin{equation}
\rho[p_{\lambda_1},p_{\lambda_2}] = \frac{p(\omega;\bar\lambda)}{\sqrt{p(\omega,\lambda_1)p(\omega,\lambda_2)}},\quad\forall \omega\in\calX
\end{equation}
where $\bar\lambda=M_\theta(\lambda_1,\lambda_2)$ is a {\em generalized quasi-arithmetic mean} induced by the ordinary-to-natural parameter conversion function $\theta(\lambda)$.
Thus our method requires only {\em partial} canonical factorization of the  densities of an exponential family to get $\theta(\lambda)$.
The formula for the Bhattacharyya  distance, Hellinger  distance, and $\alpha$-divergences follow straightforwardly:
\begin{eqnarray}
D_B[p_{\lambda_1},p_{\lambda_2}] &=& \log\left( \frac{\sqrt{p(\omega,\lambda_1)p(\omega,\lambda_2)}}{p(\omega;\bar\lambda)} \right),\\
D_H[p_{\lambda_1},p_{\lambda_2}]
&=& \sqrt{1-
\frac{p(\omega;\bar\lambda)}{\sqrt{p(\omega,\lambda_1)p(\omega,\lambda_2)}}
},\\
D_\alpha[p_{\lambda_1}:p_{\lambda_2}]  &=& \frac{1}{\alpha(1-\alpha)}  \left(1-\frac{p(\omega;\bar\lambda)}{\sqrt{p(\omega,\lambda_1)p(\omega,\lambda_2)}}\right), \quad \alpha\in\bbR\backslash\{0,1\}. 
\end{eqnarray}

In practice, it is easy to {\em program} those formula using legacy software APIs which offer many parametric densities in their library:
First, we check that the distribution is an exponential family.
Then we set $\omega$ to be {\em any} point of the support $\calX$, {\em partially factorize} the distribution in order to retrieve $\theta(\lambda)$ and its reciprocal function $\lambda(\theta)$, and equipped with these functions, we implement the corresponding {\em generalized weighted quasi-arithmetic   mean} function $M_{\theta,\alpha}(a,b)=\theta^{-1}(\alpha\theta(a)+(1-\alpha)\theta(b))$ to calculate $\bar\lambda=M_\theta(\lambda_1,\lambda_2)$.  

To calculate the  Kullback-Leibler divergence (and Jeffreys' divergence) without the explicit use of the cumulant function, we reported two methods:
The first method consists in expressing the KLD as a limit of $\alpha$-skew Bhattacharyya distance which writes as:
\begin{equation}
D_\KL[p_{\lambda_1}:p_{\lambda_2}] =  \log\left(\frac{p(\omega;\lambda_1)}{p(\omega;\lambda_2)}\right) 
+ \lim_{\alpha\rightarrow 0}   \frac{1}{\alpha} \log\left( \frac{p\left(\omega;M_{\theta,\alpha}(\lambda_2,\lambda_1)\right)}{p(\omega,\lambda_1)}\right).
\end{equation}
This limit can be calculated symbolically using a computer algebra system, or approximated for a small value of $\alpha$ by
\begin{equation}
D_{\KL,\alpha}[p_{\lambda_1}:p_{\lambda_2}] = \log\left(\frac{p(\omega;\lambda_1)}{p(\omega;\lambda_2)}\right) 
+   \frac{1}{\alpha} 
\log\left(\frac{
p\left(\omega;M_{\theta,\alpha}(\lambda_2,\lambda_1)\right)}{p(\omega;\lambda_1)}\right).
\end{equation}

When dealing with uni-order exponential family, we can use a first-order approximation of the weighted quasi-arithmetic mean to express the KLD as the following limit:
\begin{equation}
D_\KL[p_{\lambda_1}:p_{\lambda_2}] = \log\left(\frac{p(\omega;\lambda_1)}{p(\omega;\lambda_2)}\right)
+ \lim_{\alpha\rightarrow 0} \frac{1}{\alpha} 
\log\left(\frac{
p\left(\omega;\lambda_1+\alpha \frac{\theta(\lambda_2)-\theta(\lambda_1)
}{\theta'(\lambda_1)}\right)}{p(\omega;\lambda_1)}\right).
\end{equation}

Notice that we can also estimate $D_{\KL,\alpha}$, $\rho_\alpha$  and related dissimilarities (e.g., when the cumulant function is intractable) using density ratio estimation techniques~\cite{densityrationestimation-2012}.

The second approach consists in using the entropy and moment formula which are often available when dealing with parametric distributions.
When the parametric distributions form an exponential family, the KLD is equivalent to a Legendre-Fenchel divergence, and we write this  Legendre-Fenchel divergence as:
\begin{equation}
D_\KL[p_{\lambda_1}:p_{\lambda_2}] = -\log p(\omega;\lambda_2) -h(p_{\lambda_1}) -E_{p_{\lambda_1}}[k(x)] -\theta(\lambda_2)^\top E_{p_{\lambda_1}}[t(x)].
\end{equation}
It follows that the Jeffreys' divergence is expressed as
\begin{equation}
D_J[p_{\lambda_1},p_{\lambda_2}] = \left(\theta(\lambda_2)-\theta(\lambda_1)\right)^\top (E_{p_{\lambda_2}}[t(x)] - E_{p_{\lambda_1}}[t(x)]).
\end{equation}

Finally, we proved in \S\ref{sec:KLdr} that the Kullback-Leibler divergence between two densities $p_{\lambda_1}$ and $p_{\lambda_2}$ of an exponential family $\calE$ of order $D$ can be expressed as
 $$
D_\KL[p_{\lambda_1}:p_{\lambda_2}] = \frac{1}{s}\sum_{i=1}^s \log\left( \frac{p_{\lambda_1}(\omega_i)}{p_{\lambda_2}(\omega_i)} \right),
$$
 where
$\omega_1,\ldots, \omega_s$ are $s\leq D+1$ distinct samples of $\calX$ chosen such that 
$\frac{1}{s}\sum_{i=1}^s t(\omega_i)=E_{p_{\lambda_1}}[t(x)]$.
We illustrated how to find the $\omega_i$'s for the univariate Gaussian family and the multivariate zero-centered Gaussian family.

To conclude this work, let us emphasize that we have revealed a new kind of {\em invariance} when providing closed-form formula for
 common (dis)similarities between densities of an exponential family without explicitly using the cumulant function of that exponential family:
For the Bhattacharrya/Hellinger/$\alpha$-divergences, the $\omega$ can be chosen as {\em any} arbitrary point of the support $\calX$. For the Kullback-Leibler divergence, by carefully choosing a set of $\omega$'s, we may express the Kullback-Leibler divergence as a
 weighted sum of log density ratios.

\bibliographystyle{plain}
\bibliography{BhatTrickBIBV3}

\appendix

\section*{Closed-form formula using the {\sc Maxima} computer algebra system}\label{sec:CAS}

Since the statistical (dis)similarities rely on integral calculations, we may use symbolic calculations to check the results.
For example, below is some code snippets written in {\sc Maxima}.\footnote{http://maxima.sourceforge.net/}
The code snippet below calculates symbolically the Bhattacharyya coefficient for several exponential families.

\begin{verbatim}
/* Quasi-arithmetic mean associated with the univariate Gaussian family */
ptheta(lambda):=[lambda[1]/lambda[2],-1/(2*lambda[2])];
plambda(theta):=[-theta[1]/(2*theta[2]),-1/(2*theta[2])];
ptheta(plambda([t0,t1]));
l1: [p1,p2];
l2: [q1,q2];
plambda(0.5*ptheta(l1)+0.5*ptheta(l2));
ratsimp(%);
/* end */

/* Quasi-arithmetic mean associated with the inverse Gaussian family */
ptheta(lambda):=[-lambda[2]/(2*lambda[1]**2),-lambda[2]/2];
plambda(theta):=[sqrt(theta[2]/theta[1]),-2*theta[2]];
ptheta(plambda([t0,t1]));
l1: [p1,p2];
l2: [q1,q2];
plambda(0.5*ptheta(l1)+0.5*ptheta(l2));
ratsimp(%);
/* end */

/* Exponential family of exponential distributions */
assume(lambda1>0);
assume(lambda2>0);
p(x,lambda) := lambda*exp(-lambda*x); 

integrate(sqrt(p(x,lambda1)*p(x,lambda2)),x,0,inf);
ratsimp(%);
/* end */

/* Exponential family of zero-centered Gaussian densities */
assume(sigma>0);
p(x,sigma) :=  (1.0/(2*sigma))*exp(-abs(x)/sigma);
assume(sigma1>0);
assume(sigma2>0);

integrate(sqrt(p(x,sigma1)*p(x,sigma2)),x,-inf,inf);
ratsimp(%);
/* end */

/* Exponential family of centered-Laplacian distributions */
assume(lambda1>0);
assume(lambda2>0);
p(x,lambda) := (1/(2*lambda))*exp(-abs(x)/lambda);

integrate(sqrt(p(x,lambda1)*p(x,lambda2)),x,-inf,inf);
ratsimp(%);
/* end*/

/* Exponential family of Weibull densities with prescribed shape parameter k */
declare( k , integer);
assume(k>=1);
assume(lambda1>0);
assume(lambda2>0);
p(x,lambda) := (k/lambda)*(x/lambda)**(k-1)*exp(-(x/lambda)**k);

integrate(sqrt(p(x,lambda1)*p(x,lambda2)),x,0,inf);
expand(ratsimp(%));
/* end */

/* KLD betwen Weibull distributions by symbolic computing of the limit */
declare( k , integer);
assume(lambda1>0);
assume(lambda2>0);
k:3;
omega:1;
t(u):=u**(-k);
tp(u):=k*u**(-k-1);
p(x,l):=(k/l)*((x/l)**(k-1))*exp(-(x/l)**k);
mean(l1,l2):=l1+alpha*(t(l1)-t(l2))/tp(l1);
log(p(omega,l1)/p(omega,l2)) + (1.0/alpha)*log(p(omega,mean(l1,l2))/p(omega,l1));
limit (ratsimp(%), alpha, 0);
expand(%);
/* end */
\end{verbatim}

\section{Further illustrating examples}\label{sec:furtherexamples}

The Laplacian  exponential family illustrates the use of the harmonic mean $H(a,b)=\frac{2ab}{a+b}$ for $a,b>0$:

\begin{example}\label{ex:Laplacianfamily}
Consider the family of zero-centered Laplacian densities~\cite{EF-2009} $p_\lambda(x)=\frac{1}{2\lambda}\exp(-\frac{|x|}{\lambda})$ 
for $x\in\bbR$.
We have $p_\lambda(\omega)=\frac{1}{2\lambda}$ for $\omega=0$, $\theta(u)=-\frac{1}{u}=\theta^{-1}(u)$ so that
$\bar\lambda=M_\theta(\lambda_1,\lambda_2)=H(\lambda_1,\lambda_2)=\frac{2\lambda_1\lambda_2}{\lambda_1+\lambda_2}$, where $H$ denotes the {\em harmonic mean}.
Applying the cumulant-free formula Eq.~\ref{eq:bhatcoeffcf}, we get
\begin{eqnarray*}
\rho[p_{\lambda_1},p_{\lambda_2}] &=& \frac{\sqrt{p(\omega,\lambda_1)p(\omega,\lambda_2)}}{p\left(\omega;\bar\lambda\right)},\\
&=& \frac{1}{2 \frac{2\lambda_1\lambda_2}{\lambda_1+\lambda_2}} \frac{1}{\sqrt{\frac{1}{2\lambda_1}\frac{1}{2}\lambda_2}},\\
&=& \frac{2\sqrt{\lambda_1\lambda_2}}{\lambda_1+\lambda_2}.
\end{eqnarray*}
Note that the arithmetic mean $A(\lambda_1,\lambda_2)=\frac{\lambda_1+\lambda_2}{2}$ dominates the geometric mean $G(\lambda_1,\lambda_2)=\sqrt{\lambda_1\lambda_2}$ (i.e., $A\geq G$) so that $\rho[p_{\lambda_1},p_{\lambda_2}]= \frac{G(\lambda_1,\lambda_2)}{A(\lambda_1,\lambda_2)}\in (\omega,1]$.
It follows that the Bhattacharyya distance is
$$
D_\Bhat[p_{\lambda_1},p_{\lambda_2}] = \log \frac{\lambda_1+\lambda_2}{2} - \log\sqrt{\lambda_1\log \lambda_2}.
$$
This yields the same formula as the exponential distribution.
\end{example}

\begin{example}
The univariate Gaussian distribution with source parameters $\lambda=(\mu,\sigma^2)$ 
has density $p_\lambda(x)=\frac{1}{\sqrt{2 \pi \lambda_2}} \exp \left(-\frac{(x-\lambda_1)^{2}}{2 \lambda_2}\right)$.
The have $\theta(u)=(\frac{u_1}{u_2},-\frac{1}{2u_2})$ and $\theta^{-1}(u)=(-\frac{u_1}{2u_2},-\frac{1}{2u_2})$.
We have
\begin{equation}
M_\theta(\lambda_1,\lambda_2)=\left(\frac{\mu_1\sigma_2^2+\mu_2\sigma_1^2}{\sigma_1^2+\sigma_2^2},\frac{2\sigma_1^2\sigma_2^2}{\sigma_1^2+\sigma_2^2}\right).
\end{equation}
\end{example}

Our next example illustrates an unusual non-separable mean derived from the inverse Gaussian family:

\begin{example}\label{ex:InverseGaussianBhat}
Consider the family of inverse Gaussian densities:
$$
p(x ; \mu, \lambda)=\sqrt{\frac{\lambda}{2 \pi x^{3}}} \exp \left(-\frac{\lambda(x-\mu)^{2}}{2 \mu^{2} x}\right),
$$
with $\calX=(0,\infty)$, and source parameters $(\mu>0,\lambda>0)$.
The inverse Gaussian densities form an exponential family with natural parameters $\theta(\lambda)=\left(-\frac{\lambda_2}{2\lambda_1^2},-\frac{\lambda_2}{2}\right)$.
The inverse conversion function is $\theta^{-1}(\theta_1,\theta_2)=\left(\frac{\theta_2}{\theta_1},-2\theta_2\right)$.
Thus the generalized quasi-arithmetic mean induced by the multivariate $\theta(\cdot)$ is 
$$
M_\theta(\lambda,\lambda') = \left(\sqrt{ 
\frac{(\lambda_2+\lambda_2') \lambda_1^2(\lambda_1')^2}{\lambda_2(\lambda_1')^2+\lambda_2'\lambda_1^2} 
 } 
 ,\frac{\lambda_2+\lambda_2'}{2}  \right).
$$
The calculation can be checked using a computer algebra system (detailed in Appendix~\ref{sec:CAS}).
We choose $\omega=1$.
\end{example}

Next, we consider the case of the exponential family of central Wishart densities~\cite{WishartKMLE-2014}.

\begin{example}\label{ex:WishartBhat}
The density of a Wishart distribution~\cite{Wishart-1928} with positive-definite scale matrix $S\succ 0$ and number of degrees of freedom $n\geq d$ (matrix generalization of the chi-squared distribution) is:
\begin{equation}
p(x;\lambda=(\lambda_s=n,\lambda_M=S))=\frac{|X|^{\frac{n-d-1}{2}}
 e^{-\frac{1}{2} \operatorname{tr}\left(S^{-1} X\right)}}{2^{\frac{n d}{2}}|S|^{\frac{n}{2}} \Gamma_{d}\left(\frac{n}{2}\right)},
\end{equation}
where $\Gamma_d$ denotes the multivariate Gamma function:
\begin{equation}
\Gamma_{d}(x):=\pi^{d(d-1) / 4} \prod_{j=1}^{d} \Gamma(x+(1-j) / 2).
\end{equation}

The sample space $x\in\calX$ is the set of $d\times d$ positive-definite matrices.
The natural parameter consists of a scalar part $\theta_s$ and a matrix part $\theta_M$:
\begin{equation}
\theta(\lambda)=\left(\theta_{s}, \theta_{M}\right)=\left(\frac{\lambda_s-d-1}{2}, -\frac{1}{2}\lambda_M^{-1}\right).
\end{equation}
The inverse conversion natural$\rightarrow$ordinary function is
\begin{equation}
\lambda(\theta)=\left(2\lambda_s+d+1,-\frac{1}{2}\lambda_M^{-1}\right).
\end{equation}

It follows that the quasi-arithmetic mean is the following separable scalar-vector quasi-arithmetic mean (scalar arithmetic mean with a harmonic matrix mean):
\begin{equation}
M_\theta((n_1,S_1),(n_2,S_2))=\left(\frac{n_1+n_2}{2},\left(\frac{S_1^{-1}+S_2^{-1}}{2}\right)^{-1} \right).
\end{equation}

We choose $\omega=I$ (the identity matrix) so that 
\begin{equation}
p(I;S,n) =\frac{ 
 \exp\left( -\frac{1}{2} \tr(S^{-1}) \right) 
}{ 
2^{\frac{n d}{2}}  |S|^{\frac{n}{2}} \Gamma_{d}\left(\frac{n}{2}\right)
}.
\end{equation}

The sufficient statistics is $t(x)=(\log|X|,-\frac{1}{2}X)$.
It follows that $\eta=\nabla F(\theta)=E[t(x)]=(-\log |\frac{S^{-1}}{2}|+\sum_{i=1}^d \psi(\frac{n+1-i}{2}), -\frac{1}{2}nS)$ 
(see~\cite{VariationalBayes-2006}),
where $\psi_d$ is the multivariate digamma function (the derivative of the logarithm of the multivariate gamma function).
The differential entropy is: 
$$
h(p_{n,S})=\frac{d+1}{2} \ln |S|+\frac{1}{2} d(d+1) \ln 2+\ln \Gamma_{d}\left(\frac{n}{2}\right)-\frac{n-d-1}{2} \psi_{d}\left(\frac{n}{2}\right)+\frac{n d}{2}.
$$

It follows that the Kullback-Leibler divergence between two central Wishart densities is:

\begin{equation}
D_\KL\left(p_{n_1,S_1}: p_{n_2,S_2}\right)=-\log \left(\frac{\Gamma_{d}\left(\frac{n_{1}}{2}\right)}{\Gamma_{d}\left(\frac{n_{2}}{2}\right)}\right)+\left(\frac{n_{1}-n_{2}}{2}\right) \psi_{d}\left(\frac{n_{1}}{2}\right)+
\frac{n_{1}}{2}\left(-\log \frac{\left|S_{1}\right|}{\left|S_{2}\right|}+\operatorname{tr}\left(S_{2}^{-1} S_{1}\right)-d\right).
\end{equation}

\end{example}

This last example exhibits a non-usually mean when dealing with the Bernoulli family:

\begin{example}\label{ex:Bernoullifamily}
Consider the Bernoulli family with density $p_\lambda(x)=\lambda^x(1-\lambda)^{1-x}$ where $\calX=\{0,1\}$ and 
$\lambda$ denotes the probability of success. The Bernoulli family is a discrete exponential family.
We have $\theta(u)=\log\frac{u}{1-u}$ and $\theta^{-1}(u)=\frac{e^u}{1+e^u}$.
The associated quasi-arithmetic mean is 
\begin{equation}
M_\theta(\lambda_1,\lambda_2) = \sqrt{ \frac{\lambda_1\lambda_2}{(1-\lambda_1)(1-\lambda_2)}}.
\end{equation}
\end{example}

\end{document}